\documentclass[a4paper,11pt]{amsart}
\usepackage{amscd,amsmath,amsthm}
\newtheorem{theorem}{Theorem}[section]
\newtheorem{cor}[theorem]{Corollary}
\newtheorem{lemma}[theorem]{Lemma}
\newtheorem{prop}[theorem]{Proposition}
\theoremstyle{remark}
\newtheorem{remark}[theorem]{Remark}
\newtheorem{remarks}[theorem]{Remarks}
\newtheorem{example}[theorem]{Example}
\newtheorem{examples}[theorem]{Examples}
\theoremstyle{definition}
\newtheorem{definition}[theorem]{Definition}
\numberwithin{equation}{section}
\DeclareMathOperator{\Hom}{Hom}
\DeclareMathOperator{\cv}{cov}

\DeclareMathOperator{\id}{id}

\DeclareMathOperator{\Span}{span}
\DeclareMathOperator{\clsp}{\overline{span}}
\newsymbol\rtimes 216F
\newcommand{\abs}[1]{\lvert#1\rvert}
\newcommand{\norm}[1]{\lVert#1\rVert}
\newcommand{\cstar}{$C^*$\ndash}
\newcommand{\cov}{\text{$C^*_{{\cv}}$}}
\newcommand{\Star}{${}^*$\ndash}
\newcommand{\cross}{\rtimes}
\newcommand{\bpp}{B_P \textstyle{\cross_{\tau, E}} P}
\newcommand{\ap}{A \textstyle{\cross_{\beta, E}} P}
\newcommand{\rankone}[2]{#1 \otimes\overline{#2}}
\newcommand{\set}[1]{\{#1\}}
\newcommand{\ip}[2]{\langle #1, #2 \rangle}
\newcommand{\bh}{\Bb(\Hh)}
\newcommand{\ra}{\rightarrow}

\newcommand{\ndash}{\nobreakdash-}
\newcommand{\Ndash}{\nobreakdash--}
\newcommand{\field}[1]{\mathbb{#1}}
\newcommand{\CC}{\field{C}}
\newcommand{\NN}{\field{N}}

\newcommand{\RR}{\field{R}}
\newcommand{\TT}{\field{T}}
\newcommand{\ZZ}{\field{Z}}

\newcommand{\Bb}{{\mathcal B}}
\newcommand{\Cc}{{\mathcal C}}
\newcommand{\Ff}{{\mathcal F}}
\newcommand{\Gg}{{\mathcal G}}
\newcommand{\Hh}{{\mathcal H}}
\newcommand{\Kk}{{\mathcal K}}
\newcommand{\Oo}{{\mathcal O}}
\newcommand{\Pp}{{\mathcal P}}
\newcommand{\Tt}{{\mathcal T}}
\newcommand{\Uu}{{\mathcal U}}
\begin{document}
%
%
\title[discrete product systems]{Discrete Product Systems and\\ Twisted Crossed
Products by Semigroups}
\author[Neal Fowler]{Neal Fowler}
\address{Department of Mathematics  \\
      University of Newcastle\\  NSW  2308\\ AUSTRALIA}
\email{neal@math.newcastle.edu.au}
\author[Iain Raeburn]{Iain Raeburn}
\email{iain@math.newcastle.edu.au}
\date{\today}
\thanks{This research was supported by the Australian Research 
Council. It was completed while the second author was visiting Dartmouth College;
he thanks Dana Williams and his colleagues for their warm hospitality and their
excellent facilities.} 
\begin{abstract}
A product system $E$ over a semigroup $P$ is a family of Hilbert
spaces $\{E_s:s\in P\}$ together with multiplications $E_s\times E_t\to
E_{st}$. We view $E$ as a unitary-valued cocycle on $P$, and consider
twisted crossed products $\ap$ involving $E$ and an
action $\beta$ of $P$ by endomorphisms of a $C^*$-algebra $A$. When $P$
is quasi-lattice ordered in the sense of Nica, we isolate a class of
covariant representations of $E$, and consider a  twisted crossed
product $\bpp$ which is universal for covariant
representations of $E$ when $E$ has finite-dimensional fibres, and in
general is slightly larger. In particular, when $P=\NN$ and $\dim
E_1=\infty$, our algebra $B_\NN\rtimes_{\tau,E}\NN$ is a new infinite
analogue of the Toeplitz-Cuntz algebras ${\mathcal TO}_n$. Our main theorem
is a characterisation of the faithful representations of $\bpp$.
\end{abstract}
\maketitle
%
%
Crossed  products of \cstar algebras by semigroups of endomorphisms 
have been profitably used to model Toeplitz algebras
\cite{alnr, adji, lacarae} and the Hecke algebras arising in the
Bost-Connes analysis of phase transitions in number theory
\cite{lacarae2, alr, laca2}. There are two main ways of studying such a
crossed product. First, one can try to embed it as a corner in a
crossed product by an automorphic action of an enveloping group, and
then apply the established theory.
The algebra on which the group acts is typically a direct limit, and the
success of this approach depends on being able to recognise the direct
limit and the action on it \cite{cuntz, sta, murnew}. Or, second, one
can use the techniques developed in \cite{bkr, alnr, lacarae} to deal
directly with the semigroup crossed product and its representation
theory. Here the goal is a characterisation of the faithful
representations of the crossed product, and such characterisations have
given important information about a wide range of semigroup
crossed products \cite{alnr, lacarae, lacarae2, alr}.  

For ordinary crossed products $A\rtimes_\alpha G$ (those
involving an action $\alpha$ of $G$ by automorphisms of $A$), an important
adjunct are the twisted crossed products
$A\rtimes_{\alpha,\omega} G$, in which the multiplication of elements of
$G$ has been twisted by a cocycle $\omega$. This cocycle might take
values in the unitary groups of $A$, $M(A)$ or $ZM(A)$, but the most
important are the scalar-valued cocycles $\omega:G\times G\to \TT$.
There is no obvious technical obstruction to developing a theory of
twisted semigroup crossed products, and indeed this has already been
done by Laca for scalar-valued cocycles on totally ordered groups
\cite{laca1}. Since scalar-valued cocycles on semigroups often extend to
the enveloping group \cite{lacarae0}, one might expect this theory to
be a routine combination of ideas involving semigroup crossed products
and ordinary twisted crossed products. 

In this paper we investigate a  phenomenon which arises
only for semigroups: crossed products twisted by unitary
cocycles acting on Hilbert spaces of varying dimension. Such cocycles
were introduced by Arveson under the name of product systems \cite{arv}.
The idea is to associate to each element $s$ of the semigroup $S$ a Hilbert space
$E_s$, and then the cocycle describes a multiplication from 
$E_s\times E_t$ to $E_{st}$; a scalar cocycle
$\omega:S\times S\to \TT$ determines such a system by
taking $E_s=\CC$ for all $s$ and using $(w,z)\mapsto
\omega(s,t)wz$ as the product from $E_s\times E_t$ to $E_{st}$.
Because $\dim E_{st}=\dim E_s\times \dim E_t$, product systems with
fibres of  dimension other than $1$ cannot exist on groups (at least
in a naive sense), so the possibility of twisting crossed products by
product systems is appropriate only for actions of semigroups. It is not
an entirely new idea: the crossed products of multiplicity $n$ of
Stacey \cite{sta} are twisted crossed products by actions of the
semigroup $\NN$ in which the product system $E$ has $\dim E_1=n$.

Of the various kinds of semigroups studied in the literature, we have
chosen to work with the quasi-lattice ordered semigroups of Nica
\cite{nica}; these include the totally ordered groups considered in
\cite{murjot, alnr, laca1}, the direct sums $\NN^k$, and the free
products considered in \cite{lacarae}. For a product system $E$ over
such a semigroup $P$, one can define a natural notion of covariant
representation generalising that of \cite{nica, lacarae}: loosely
speaking, a representation $\phi$ of $E$ is a family of isometric maps
$\phi_s:E_s\to B(\Hh)$ such that each $\phi_s(v)$ is an isometry
and $\phi_{st}(uv)$ is the composition of the operators $\phi_s(u)$ and
$\phi_t(v)$, and $\phi$ is covariant if the projections on the ranges
$\phi_s(E_s)$ are aligned in a manner compatible with the ordering on
$P$. The motivating example is the trivial product system on $\NN^2$,
where the representations are given by two commuting isometries and the
covariant representations by two \Star commuting isometries.

The main results of \cite{nica} and \cite{lacarae} concern the
\cstar algebra, here denoted $\cov(P)$, which is universal for
covariant isometric representations of the quasi-lattice ordered group
$(G,P)$. In \cite{lacarae}, $\cov(P)$ is viewed as a semigroup
crossed product $B_P\rtimes_\tau P$, where $B_P$ is the
\cstar subalgebra of $\ell^\infty(P)$ spanned by the characteristic
functions $1_x:=\chi_{xP}$, and $\tau_t(1_x)=1_{tx}$. Here we aim to view
the universal \cstar algebra $\cov(P,E)$ for covariant
representations of $E$ as a twisted crossed product
$\bpp$, and use techniques like those of \cite{lacarae}
to characterise their faithful representations. However, carrying out
this program has raised some intriguing issues. 

We shall construct suitable twisted crossed products
$\bpp$, and show that the \cstar subalgebra of
$\bpp$ generated by the canonical copy  of $E$ is
universal for covariant representations of $E$, and hence can
reasonably be denoted $\cov(P,E)$. When the fibres of $E$ are
finite-dimensional, $\cov(P,E)$ is all of $\bpp$, but
in general it may not be. This last phenomenon occurs, for example, when
$P=\NN$ and $E_1$ is infinite-dimensional: $\cov(\NN,E)$ is the
Cuntz algebra $\Oo_\infty$ generated by isometries
$\{V_k:k\in\NN\}$ with orthogonal ranges, whereas $B_\NN\rtimes_{\tau,E}\NN$
contains the projection $1-\sum_{k=1}^\infty V_kV_k^*=1-\tau_1(1)$.
This undermines the popular view that the Cuntz algebra
$\Oo_\infty$ coincides with the Toeplitz-Cuntz algebra $\Tt\Oo_\infty$,
since
$B_\NN\rtimes_{\tau,E}\NN$ seems a logical candidate for the latter.
Our main theorem characterises faithful representations of
$\bpp$ rather than $\cov(P,E)$, and thus achieves
our goal only for systems with finite-dimensional fibres. We plan to
return to the topic of systems with infinite-dimensional fibres in a
sequel.

We have organised our work as follows. We begin with an introductory
section on product systems and their representations, giving a variety
of examples and constructions. General twisted crossed products are
discussed only in \S2: as in \cite{lacarae}, we are mainy interested in
the specific crossed products $\bpp$ which capture the
covariance condition on representations of $E$. The covariance condition
itself is modelled on that of Nica, and only makes sense for product
systems on quasi-lattice ordered semigroups. In \S3 we discuss it and
its connection with covariant representations of the system
$(B_P,P,\tau,E)$. We can then prove that $\cov(P,E)$ embeds
naturally in the semigroup crossed product $\bpp$
(Theorem~\ref{theorem:subalgebra}). 

Our main theorem is our characterisation of faithful representations of
$\bpp$. There are two main steps. First, under an
amenability hypothesis, we follow the procedure pioneered by Cuntz,
which reduces the problem to proving an estimate concerning the
deletion of off-diagonal terms. The details are necessarily different,
but the general plan of \cite[\S3]{lacarae} carries over under a
spanning hypothesis on the product system which holds in the
interesting examples. Second, we have to verify the amenability
hypothesis in a reasonable number of situations. It is automatic, for
example, if the enveloping group of $P$ is amenable, or if $P$ is a
free product of such semigroups and the product system satisfies a
modest-looking spanning condition. Both the spanning conditions we have
mentioned are satisfied if $E$ has finite-dimensional fibres, so our
main theorem applies to all such product systems on $\NN^k$ or on
free products of subsemigroups of amenable groups.

\section{Product Systems and their Representations}\label{section:ps}

\begin{definition}\label{defn:ps}
Suppose $P$ is a semigroup with identity and $p:E\to P$ is a family of
nontrivial complex Hilbert spaces whose fibre over the identity
is one-dimensional.  Write $E_t$ for the fibre
$p^{-1}(t)$ over $t\in P$.  We say that
$E$ is a {\em (discrete) product system over $P$\/} if $E$ is a semigroup, $p$ is a
semigroup homomorphism, and for each $s,t\in P$ the map
$(u,v) \in E_s \times E_t \mapsto uv \in E_{st}$
extends to a unitary isomorphism $U_{s,t}$ of $E_s\otimes E_t$ onto $E_{st}$.
\end{definition}

\begin{remark}
The associativity of multiplication in the semigroup $E$ implies that the unitary
operators
$U_{s,t}$ satisfy
\[
U_{rs, t}(U_{r,s} \otimes I) = U_{r,st}(I \otimes U_{s,t})
\]
for $r,s,t\in P$.
Thus  product systems over $P$
can be viewed as unitary $2$-cocycles acting on a varying but coherent system of 
Hilbert spaces.
\end{remark}

\begin{lemma}\label{lemma:unital} Suppose $E$ is a product system over
a semigroup $P$ with identity $e$.  Then $E$
has an identity $\Omega$ such that
$p(\Omega) = e$ and $\norm\Omega = 1$.
\end{lemma}

\begin{proof} Let $z$ be a unit vector in $E_e$.
Then $z^2 \in E_e$ also, so $z^2 = \lambda z$ for some $\lambda\in\CC$ such that
$\abs\lambda^2 = \ip{\lambda z}{\lambda z}
  = \ip{z^2}{z^2} = \ip zz \, \ip zz = 1$.
Suppose $x\in E$.  Then $zx\in E_{p(x)}$, and for any $y\in E_{p(x)}$,
\[
\ip{zx}y = \ip zz\ip{zx}y
  = \ip{(zz)x}{zy} = \ip{\lambda zx}{zy}
  = \ip zz \ip{\lambda x}y = \ip{\lambda x}y,
\]
so $zx = \lambda x$.  Similarly we have $xz = \lambda x$,
and thus $\Omega = \overline\lambda z$
is an identity for $E$.  We have  $p(\Omega) = e$ because $z\in E_e$, 
and $\norm\Omega = 1$ because $z$ is a unit vector.
\end{proof}

\begin{examples}\label{examples:ps}
(E1) The {\em trivial product system\/} over $P$
is the trivial bundle $P \times \CC$ with multiplication given by
$(s, w)(t, z) = (st, wz)$.

(E2) ({\em Lexicographic Product Systems\/})
Given a product system $p:E\to P$ with $\dim E_t < \infty$
for each $t\in P$,
the dimension function $d: t\mapsto \dim E_t$ is a semigroup homomorphism of $P$
into the multiplicative positive integers $\NN^*$.
Conversely, given $d\in\Hom(P,\NN^*)$,
we can construct a product system over $P$
with dimension function  $d$ as follows.
Let $E = \bigsqcup_{t\in P} \set{t} \times \CC^{d(t)}$,
$p(t,v) = t$, and define multiplication in $E$
by $(s,u)(t,v) = (st, w)$ where
\[
w_{(i-1)d(t) + j} = u_iv_j, \qquad 1 \le i \le d(s),\ 1 \le j \le d(t).
\]
Since this construction is based on the lexicographic ordering of
\[
\set{1,2,\dots, d(s)} \times \set{1,2,\dots, d(t)},
\]
we call $E$ the {\em lexicographic product system over $P$ determined by $d$}.

(E3) Suppose $p:E\to P$ is a
product system over a semigroup $P$ with identity $e$,
and $\mu$ is a multiplier on $P$;
that is, $\mu:P\times P \to \TT$ satisfies
\begin{itemize}
\item $\mu(t,e) = 1 = \mu(e,t)$ for each $t\in P$, and
\item $\mu(r,s)\mu(rs,t) = \mu(s,t)\mu(r,st)$ for each $r,s,t\in P$;
\end{itemize}
alternatively, one might say $\mu$ is a $2$-cocycle on $P$ with values in $\TT$.
Let $E^\mu = E$, $p^\mu = p$, and define
multiplication by $(u,v) \mapsto \mu(p(u), p(v)) uv$. Then $E^\mu$ is a product
system over $P$; we say that $E^\mu$ is {\em $E$ twisted by $\mu$\/}.

If $\nu$ is another multiplier on $P$, then  $E^\mu$
is isomorphic to  $E^\nu$ iff $[\mu] = [\nu]$
as elements of the second cohomology group $H^2(P, \TT)$; the automorphism group of
$E^\mu$ is $\Hom(P, \TT)$.

(E4) 
For each $\lambda$, let $p:E^\lambda\to P^\lambda$ be a product system. Then there is
a product system $*E^\lambda$ over the free product $*P^\lambda$: for a reduced word
$s=s_1\dotsm s_n\in*P^\lambda$, say $s_i\in P^{\lambda_i}$, we take
$(*E^\lambda)_s:=E^{\lambda_1}_{s_1}\otimes\dotsb\otimes E^{\lambda_n}_{s_n}$, and
if also $t=t_1\dotsm t_m\in*P^\lambda$, say $t_i\in P^{\mu_i}$, we define
\[
(w_1\otimes\dotsb\otimes w_n)(v_1\otimes\dotsb\otimes v_m):=
\begin{cases}
w_1\otimes\dotsb\otimes w_nv_1\otimes\dotsb\otimes v_m&\text{if $\lambda_n=\mu_1$}\\
w_1\otimes\dotsb\otimes w_n\otimes v_1\otimes\dotsb\otimes v_m&\text{otherwise}.
\end{cases}
\]
\end{examples}

Product systems over $\NN$ are particularly easy to describe:

\begin{prop}\label{prop:ps over N} Suppose $E$ and $F$ are product systems
over $\NN$.  Then $E$ and $F$ are isomorphic iff $E_1 \cong F_1$.
\end{prop}

\begin{proof} If $U$ is a unitary isomorphism of $E_1$ onto $F_1$, then  the unitary
operators
$U^{\otimes n} : E_1^{\otimes n} \to F_1^{\otimes n}$ induce a family of unitaries
$\psi_n: E_n \to F_n$ such that
\[
\psi_n(u_1u_2\dotsm u_n)  =  (Uu_1)(Uu_2)\dotsm (Uu_n), 
\qquad u_1, \dots, u_n \in E_1,
\]
and these combine to give an isomorphism of product systems.
\end{proof}

\begin{cor}\label{cor:ps over N} For each $d \in \set{1,2,\dots, \aleph_0}$
there is, up to isomorphism, a unique product system $E^d$ over $\NN$
whose fibre over $1$ is $d$\ndash dimensional.
\end{cor}

\begin{proof} Let $d \in \set{1,2,\dots, \aleph_0}$, fix a $d$\ndash dimensional
Hilbert space $\Hh$, let $E^d = \bigsqcup_{n=0}^\infty  \set{n} \times
\Hh^{\otimes n}$, and  define $(m,u)(n,v) := (m + n, u
\otimes v)$.
\end{proof}

\begin{definition}\label{defn:rep}
A {\em representation\/} of a product system $p:E\to P$
on a Hilbert space $\Hh$ is a map $\phi:E\to\bh$ such that
\begin{itemize}
\item[(1)] $\phi(uv) = \phi(u)\phi(v)$ for every $u,v\in E$, and
\item[(2)] $\phi(v)^*\phi(u) = \ip uv I$ whenever $p(u) = p(v)$.
\end{itemize}
\end{definition}

\begin{remarks}
(1) Condition (2) implies that every operator in the range of $\phi$
is a multiple of an isometry, and that $\phi$ is linear on the fibres of $p$; see
\cite[p.8]{arv}. It also implies that $\phi$ is
isometric, hence injective; thus each vector space
$\phi(E_t)$ has a Hilbert space structure in which
the inner product is given by
$\ip ST I = T^*S$, and the corresponding Hilbert space norm on
$\phi(E_t)$ agrees with the operator norm.

(2) Condition (1) implies that $\phi(\Omega)$ is an idempotent, and condition
(2)  that it is an isometry.  Thus $\phi(\Omega)$ is the identity operator
$I$.
\end{remarks}

\begin{examples}\label{examples:reps}
(1) A representation of the trivial product system  $P\times \CC$ on $\Hh$
is a homomorphism of $P$ into the semigroup of isometries on $\Hh$.
If $P\times \CC$ is twisted by a multiplier $\mu$, the representations
 are $\mu$\ndash twisted representations of $P$
by isometries.

(2) Suppose $E$ is a product system over $\NN$
with $\dim E_1 = d$.
A representation $\phi$ of $E$ will map an orthonormal basis $\{e_i\}$ for $E_1$
to a family of $d$ isometries $S_i=\phi(e_i)$ whose ranges are mutually orthogonal,
and each such family $\{S_i\}$ determines a representation of $E$.
We call $\{S_i\}$ a {\em Toeplitz-Cuntz family\/}.

(3) Let $E$ be the lexicographic product system
over $\NN\oplus\NN$ determined by the homomorphism
$d: (m,n) \in \NN\oplus\NN \mapsto 2^m3^n \in \NN^*$.  Representations of $E$
are in one-one correspondence with pairs of Toeplitz-Cuntz families
$\set{ U_1, U_2 }$, $\set{ V_1, V_2, V_3 }$ satisfying the following
commutation relations:
\begin{equation}\label{commrel}
\begin{aligned}
U_1V_1 & = V_1U_1, \\
U_1V_2 & = V_1U_2, \\
U_1V_3 & = V_2U_1,
\end{aligned}
\qquad
\begin{aligned}
U_2V_1 & = V_2U_2, \\
U_2V_2 & = V_3U_1, \\
U_2V_3 & = V_3U_2.
\end{aligned}
\end{equation}

\end{examples}

\begin{lemma}[The Left Regular Representation] \label{lemma:lrr defn}
Suppose $p:E\to P$ is a product system and $P$ is left-cancellative.
Let
$S(E) = \bigoplus_{t\in P} E_t$.
Then there is a unique representation $l:E\to\Bb(S(E))$
such that $l(v)w = vw$ for $v,w\in E$.
\end{lemma}

\begin{proof} Suppose $v\in E$ and $w = \bigoplus w_t \in S(E)$.
Since $p(vw_s) = p(vw_t)$ only when $s = t$,
the infinite series $\sum_{t\in P} vw_t$ converges in norm to
a vector $l(v)w$ of norm $\norm v\,\norm w$,
thus defining a bounded linear operator
$l(v)$ on $S(E)$.  The associativity of the product in $E$ implies that
$l:E\to\Bb(S(E))$ is  multiplicative, and if $u,v\in E$ have $p(u) = p(v)$,
then for any $w,z\in E$
\[
\begin{split}
\ip{l(v)^*l(u)w}{z} = \ip{uw}{vz}
& = \begin{cases}
            \ip uv \ip wz & \text{if $p(w) = p(z)$} \\
            0             & \text{otherwise}
      \end{cases} \\
& = \bigl\langle \ip uv w, z \bigr\rangle,
\end{split}
\]
so that $l(v)^*l(u) = \ip uv I$.
Thus $l$ is a representation of $E$.
\end{proof}

The following two propositions introduce concepts and notation which
will be used throughout the remainder of this paper.
The first is merely a translation of (\cite{arv}, Proposition~2.7) to our setting,
so we omit the proof.

\begin{prop}\label{prop:alphaphi}
Suppose $E$ is a product system over $P$ and
$\phi:E\to\bh$ is a representation.
For each $t\in P$ there is a unique normal \Star en\-do\-mor\-phism $\alpha^\phi_t$
of $\bh$ such that
\[
\phi(E_t) = \set{ T\in\bh: \alpha^\phi_t(A)T = TA\quad\text{for each $A\in\bh$} };
\]
the map $t\mapsto \alpha^\phi_t$ is a semigroup homomorphism.
If $\Bb$ is an orthonormal basis for $E_t$, then $\alpha_t^\phi$ is given by the
strongly convergent sum
\[
\alpha^\phi_t(A) = \sum_{u\in \Bb} \phi(u)A\phi(u)^*.
\]
\end{prop}

\begin{prop}\label{prop:rho}
Suppose $E$ is a product system over $P$ and
$\phi:E\to\bh$ is a representation.

\textup{(1)} For each $t\in P$ there is a unique faithful normal \Star
homomorphism
$\rho^\phi_t: \Bb(E_t) \to \bh$ such that
\[
\rho^\phi_t(\rankone uv) = \phi(u)\phi(v)^*\ \text{ for } u,v\in E_t,
\]
where $\rankone uv$ denotes the rank-one operator
$w \mapsto \ip wv u$ on $E_t$.

\textup{(2)} If $Q$ is a nonzero projection on $\Hh$ and $t\in P$, then the map
$T \mapsto \alpha^\phi_t(Q)\rho^\phi_t(T)$is a faithful normal \Star
homomorphism.
\end{prop}

\begin{proof} (1) Let $\Bb$ be an orthonormal basis for $E_t$.
Since $\set{ \rankone uv: u,v\in\Bb }$ is a self-adjoint
system of matrix units which generate $\Bb(E_t)$ and
$\set{ \phi(u)\phi(v)^*:u,v\in\Bb }$ is also a self-adjoint
system of nonzero matrix units, the map $\rankone uv \mapsto \phi(u)\phi(v)^*$
extends to the desired homomorphism $\rho^\phi_t$.

(2) For any $u,v\in\Bb$, note that
$\alpha^\phi_t(Q)\phi(u)\phi(v)^* = \phi(u)Q\phi(v)^*$.
Since $\set{ \phi(u)Q\phi(v)^*: u,v\in\Bb }$ is a self-adjoint system
of nonzero matrix units, the map $\rankone uv \mapsto \phi(u)Q\phi(v)^*$
extends as claimed.
\end{proof}

\section{Twisted Semigroup Crossed Products}\label{section:crossed products}

In this section we discuss how to twist semigroup crossed products
by  product systems.
We consider {\em twisted systems\/} $(A, P, \beta, E)$
in which $A$ is a unital \cstar algebra, $P$ is a semigroup with identity,
$\beta$ is an action of $P$ on $A$ by endomorphisms,
and $E$ is a product system over $P$.
We emphasise that the endomorphisms $\beta_s$ need not be unital.

\begin{definition}\label{defn:E covariance}
 A {\em covariant representation\/} of $(A,P, \beta, E)$ on a Hilbert space $\Hh$
is a pair $(\pi, \phi)$ consisting of a unital representation
$\pi:A\to\bh$ and a representation $\phi:E\to\bh$
such that $\pi\circ\beta_s = \alpha^\phi_s\circ\pi$ for $s\in P$;
by Proposition~\ref{prop:alphaphi}, this is equivalent to choosing  an orthonormal
basis  $\Bb$ for $E_s$ and asking that
\begin{equation}\label{cov}
\pi(\beta_s(a)) = \sum_{v\in\Bb} \phi(v)\pi(a)\phi(v)^* \ \text{ for }s\in P,a\in A.
\end{equation}
A {\em crossed product\/} for $(A, P, \beta, E)$ is a triple
$(B, i_A, i_E)$ consisting of a \cstar alg\-ebra $B$,
a unital \Star homomorphism $i_A:A\to B$, and a \cstar morphism $i_E:E\to B$
such that

(a) there is a faithful unital representation $\sigma$ of $B$
such that $(\sigma\circ i_A, \sigma\circ i_E)$
is a covariant representation of $(A, P, \beta, E)$;

(b) for every covariant representation $(\pi,\phi)$
of $(A, P, \beta, E)$, there is a unital
representation $\pi\times\phi$ of $B$
such that $(\pi\times\phi)\circ i_A = \pi$
and $(\pi\times\phi)\circ i_E = \phi$;

(c) the \cstar algebra $B$ is generated by
$i_A(A)\cup i_E(E)$.
\end{definition}

\begin{remark}\label{remark:crossed products}
The semigroup crossed products considered in \cite{alnr}, \cite{bkr}
and \cite{lacarae} are recovered by taking $E$ to be the trivial
product system $P\times\CC$, and twisted semigroup crossed products
$A\rtimes_{\beta,\mu}P$ involving a multiplier $\mu$  by taking $E=(P\times\CC)^\mu$,
as in Examples~\ref{examples:ps} (E3). Stacey's crossed products of multiplicity $n$
\cite{sta} are recovered by taking $E$ to be the essentially unique product system
over $\NN$ with $\dim E_1=n$ (see Corollary~\ref{cor:ps over N}).
\end{remark}

\begin{remark}
Instead of condition (a) in the definition of a crossed product, one might expect to see
something more like:

(a$'$)  for every unital representation $\sigma$ of $B$,
the pair $(\sigma\circ i_A, \sigma\circ i_E)$
is a covariant representation of $(A, P, \beta, E)$.

\noindent This condition would ensure that every unital representation of $B$
came from a covariant representation of $(A,P,\beta, E)$,
so that $B$ would be truly universal for covariant representations.
When the fibres of $E$ are finite-dimensional,
conditions (a) and (a$'$) are both equivalent to the following condition
on $(i_A, i_E)$: if $a \in A$, $s\in P$ and $\set{v_1, \dots, v_n}$
is an orthonormal basis for $E_s$, then 
\begin{equation}\label{realcov}
i_A(\beta_s(a)) = \sum_{k = 1}^n i_E(v_k)i_A(a)i_E(v_k)^*.
\end{equation}
However, if $E_s$ were infinite-dimensional,
the  sum on the right of (\ref{realcov}) would have to be infinite, and
because the isometries $i_E(v_k)$ have orthogonal range projections such sums  cannot
possibly converge in the
\cstar algebra $B$. Indeed, for systems with infinite-dimensional fibres
conditions (a) and (a$'$) need not coincide.
The following example  shows
that condition (a$'$) is too much to hope for if there is to be a crossed
product for every system with
a covariant representation.  
\end{remark}

\begin{example}\label{ex:not universal}
Consider the system $(c, \NN, \tau, E^{\aleph_0})$,
where
$\tau$ is the action of $\NN$ by translation on the algebra $c$ of convergent
sequences and $E^{\aleph_0}$ is the product system over $\NN$
with $\dim E_1=\aleph_0$ (see Corollary~\ref{cor:ps over N}).
We shall show that a crossed product $B$ for $(c,\NN, \tau, E^{\aleph_0})$ does not
satisfy condition (a$'$).

Let $\set{S_k: k\in \NN}$ be a countably-infinite
collection of isometries on a Hilbert space $\Hh$
such that $\sum S_kS_k^* = I$,
and let $\set{\delta_k: k\in \NN}$ be an orthonormal
basis for $E_1$.
The formula $\phi(\delta_k) = S_k$
extends uniquely to a representation $\phi: E \to \bh$.
Define $L: c \to \bh$
by $L(a) = \left( \lim_{k \ra \infty} a_k \right) I$.
Then $(L, \phi)$ is a
covariant representation of $(c,\NN, \tau, E^{\aleph_0})$,
and $L \times \phi(B) = C^*(\set{S_k})$.

Now let $\set{T_k: k\in \NN}$ be a family
of isometries on a Hilbert space such that
$\sum T_kT_k^* < I$.
By Cuntz's theorem the map $S_k \mapsto T_k$ extends to an isomorphism
$\pi$ of $C^*(\set{S_k})$ onto $C^*(\set{T_k})$.
Let $\sigma = \pi \circ (L \times \phi)$.  The pair
$(\sigma \circ i_c, \sigma \circ i_E) = (\pi \circ L, \pi \circ \phi)$
is not covariant since
\[
\sigma \circ i_c(\tau_1(1)) = \pi \circ L(\tau_1(1)) = \pi(I) = I,
\]
whereas 
\begin{multline*}
\alpha^{\sigma\circ i_E}_1(\sigma\circ i_c(1))
 = \alpha^{\pi \circ \phi}_1(\pi \circ L(1)) \\
 = \sum_{k=1}^\infty \pi\circ\phi(\delta_k)\pi\circ\phi(\delta_k)^*
 = \sum_{k=1}^\infty T_kT_k^*
 < I.
\end{multline*}
\end{example}

\begin{prop}\label{prop:existence of cp}
If $(A, P, \beta, E)$ has a covariant representation, then it has
a crossed product $(\ap, i_A, i_E)$
which is unique in the following sense:
if $(B, i_A', i_E')$ is another crossed product for $(A, P, \beta, E)$,
then there is an isomorphism $\theta: \ap \to B$
such that $\theta\circ i_A = i_A'$ and $\theta\circ i_E = i_E'$.

\end{prop}

\begin{proof} Say that a covariant representation $(\pi,\phi)$ is {\em cyclic\/} if the
\cstar algebra $C^*(\pi,\phi)$ generated by $\pi(A)\cup\phi(E)$ acts cyclically, i.e., has
a cyclic vector. If $(\pi,\phi)$ is any covariant representation on $\Hh$, the usual
Zorn's Lemma argument shows that $\Hh$ is the direct sum of subspaces on which
$C^*(\pi,\phi)$ acts cyclically. These subspaces are then invariant for $\pi$ and
$\phi$, and the projection $Q$ onto such a subspace commutes with $\pi(A)$ and
\Star commutes with $\phi(E)$. Since compressing by $Q$ preserves the strong operator
convergence in (\ref{cov}), the pair $(Q\pi,Q\phi)$ is covariant, and is cyclic because
$C^*(Q\pi,Q\phi)=QC^*(\pi,\phi)Q$ acts cyclically on $Q\Hh$. Thus every covariant
representation is a direct sum of cyclic representations.

Let $S$ be a set of cyclic covariant representations with
the property that every cyclic covariant representation of $(A, P, \beta, E)$
is unitarily equivalent to an element in $S$.
It can be shown that such a set $S$ exists
by fixing a Hilbert space $\Hh$ of sufficiently
large cardinality (depending on the cardinalities of $A$ and $E$)
and considering only representations on $\Hh$.
Note that $S$ is nonempty because the system has a covariant representation, which has
a cyclic summand. Define $i_A = \bigoplus_{(\pi,\phi)\in S} \pi$,
$i_E = \bigoplus_{(\pi,\phi)\in S} \phi$,
and let $\ap$ be the \cstar algebra generated by $i_A(A) \cup i_E(E)$.
Condition (a) for a crossed product is satisfied by taking $\sigma$
to be the identity representation,
condition (b) holds since every covariant representation decomposes as
a direct sum of cyclic ones,
and condition (c) was built into the definition of $\ap$.

We now prove the uniqueness.
Condition (a) allows us to realise $\ap$ and $B$
as \cstar subalgebras of $\bh$ and $\Bb(\Hh')$ in such a way that $(i_A, i_E)$
and $(i_A', i_E')$ become covariant representations of $(A, P, \beta, E)$.
Condition (b) then gives a representation
$i_A' \times i_E':\ap \to\Bb(\Hh')$
whose range is contained in $B$ because
$(i_A' \times i_E') \circ i_A = i_A'$, $(i_A' \times i_E') \circ i_E = i_E'$,
and $\ap$ is generated by $i_A(A) \cup i_E(E)$.
From (b) and (c) we see that $(i_A \times i_E) \circ (i_A' \times i_E')$
 is the identity
on $\ap$, and similarly $(i_A' \times i_E') \circ (i_A \times i_E)$
is the identity on $B$.  Hence $\theta = i_A' \times i_E'$ is the desired isomorphism.
\end{proof}

When $P$ is a subsemigroup of a group $G$, every twisted crossed product
$\ap$ carries a {\em dual coaction\/} of $G$:

\begin{prop}\label{prop:coaction}
Suppose $(A,P,\beta,E)$ is  a twisted system which has a covariant
representation. If $P$ is a subsemigroup of a group $G$, then there is an injective
coaction
\[
\delta: \ap \to (\ap) \otimes_{\min}C^*(G)
\]
such that
\[
\delta(i_A(a)) = i_A(a) \otimes 1\ \text{ and }\ 
\delta(i_E(v))= i_E(v) \otimes i_G(p(v)).
\]
If $G$ is abelian, there is a strongly continuous action $\widehat\beta$ of
$\widehat G$ on $\ap$ such that
\[
\widehat\beta_\gamma(i_A(a))=i_A(a)\ \text{ and }\ 
\widehat\beta_\gamma(i_E(v))=\gamma(p(v))i_E(v).
\]
\end{prop}

\begin{proof} Choose a faithful unital representation $\sigma$ of
$\ap$ such that $(\sigma \circ i_A, \sigma \circ i_E)$
is a covariant representation of $(A, P, \tau, E)$, and a
unitary representation $U$ of $G$ whose integrated form
$\pi_U$ is faithful on $C^*(G)$.
Then $((\sigma \circ i_A) \otimes 1, (\sigma \circ i_E) \otimes (U \circ p))$
is a covariant representation of $(A, P, \beta, E)$,
and hence there is a representation $\rho$ of $\ap$
such that
\[
\rho(i_A(a)) = \sigma \circ i_A(a) \otimes I
  = (\sigma \otimes \pi_U)(i_A(a) \otimes 1)
\]
and
\[
\rho(i_E(v)) = \sigma \circ i_E(v) \otimes U(p(v))
  = (\sigma \otimes \pi_U)(i_E(v) \otimes i_G(p(v))).
\]
Since $\sigma$ and $\pi_U$ are faithful, $\sigma \otimes \pi_U$
is faithful on $(\ap) \otimes_{\min} C^*(G)$,
and we can define $\delta := (\sigma \otimes \pi_U)^{-1} \circ \rho$.

Next let $\epsilon$ be the augmentation representation of $C^*(G)$:
$\epsilon(i_G(s)) = 1$ for all $s\in G$.  Then there is a representation
$\sigma \otimes \epsilon$ of $(\ap) \otimes_{\min} C^*(G)$
on $\Hh_\sigma = \Hh_\sigma \otimes \CC$, and  checking on generators shows
that $(\sigma \otimes \epsilon) \circ \delta = \sigma$.
Thus $\delta$ is injective.  It is also easy to check on generators
that
$(\id\otimes\delta_G)\circ\delta=(\delta\otimes\id)\circ\delta$ as homomorphisms of
$\ap$ into $(\ap) \otimes C^*(G)\otimes C^*(G)$, so $\delta$ is a coaction.

The last part follows because coactions of an abelian group $G$ are in one-to-one
correspondence with actions of $\widehat G$. Alternatively, one could use the
uniqueness of the crossed product to obtain the automorphisms $\widehat\beta_\gamma$
directly, as in \cite[Remark~3.6]{lacarae}.
\end{proof}

\section{Quasi-lattice Ordered Groups}\label{section:qlo}

Suppose $P$ is a subsemigroup of a group $G$ such that $P\cap P^{-1} = \set{e}$.
Then $s \le t$ iff $s^{-1}t \in P$ defines a  partial order on $G$ which is
left-invariant in the sense that $s \le t$ iff $rs \le rt$.
Following \cite{nica} and \cite{lacarae}, we say that
$(G,P)$ is {\em quasi-lattice ordered\/}
if every finite subset of $G$ which has an upper bound in $P$ has a
least upper bound in $P$. We shall occasionally
write $\sigma A$ for the least upper bound of a subset $A$ of $P$, and write $\sigma
A=\infty$ when $A$ has no upper bound.
Our main examples will be direct sums and free products of totally-ordered
groups of the form $(\Gamma, \Gamma^+)$, where $\Gamma$ is a countable
subgroup of $\RR$ and $\Gamma^+ = \Gamma \cap [0, \infty)$.  

\begin{remark}
We shall not  use the full strength of this definition,
so our results may be slightly more general than we have claimed.
To see why, recall from \cite{nica} that a partially
ordered group $(G,P)$ is quasi-lattice ordered if and only if

(QL1) whenever $g \in G$ has an upper bound in $P$, it has a least
upper bound in $P$, and

(QL2) whenever $s,t\in P$ have a common upper bound they have
a least common upper bound.

\noindent We  make no use of condition (QL1).
All the results in \S\ref{section:qlo}\Ndash\ref{section:system}
apply to cancellative semigroups which satisfy (QL2).
In \S\ref{section:faithfulness}
it is necessary to assume that $P$ can be embedded in a group,
but it makes no difference what the group is.  The amenability results in
\S\ref{section:amenability} can be restated
in terms of a homomorphism $\theta: P\to \Pp$  into a subsemigroup of an amenable
group. 
\end{remark}

Recall from \cite{nica} and \cite{lacarae} that a representation $V$ of $P$ by
isometries on a Hilbert space is called {\em covariant\/} if
\[
V_sV_s^*V_tV_t^* = 
\begin{cases}
  V_{s \vee t} V_{s \vee t}^* & \text{if $s \vee t < \infty$} \\
  0                           & \text{otherwise.}
\end{cases}
\]
We believe the appropriate generalisation to product systems over $P$ to be: 

\begin{definition}\label{defn:covariance}
Suppose $(G,P)$ is a quasi-lattice ordered group
and $E$ is a product system over $P$.  A representation
$\phi:E\to\bh$ is {\em covariant\/} if
\[
\alpha^\phi_s(I)\alpha^\phi_t(I) =
\begin{cases}
  \alpha^\phi_{s \vee t}(I) & \text{if $s \vee t < \infty$} \\
  0                         & \text{otherwise.}
\end{cases}
\]
\end{definition}

\begin{remark}\label{remark:covariance}
If $(G,P)$ is totally ordered, then $s \le t$ implies $\alpha^\phi_t(I) \le
\alpha^\phi_s(I)$, so every representation of $E$ is covariant.
\end{remark}

\begin{prop}\label{prop:lrr covariant}
If $(G,P)$ is a quasi-lattice ordered group and
$E$ is a product system over $P$, then the left regular representation $l$
of $E$ is covariant.
\end{prop}

For the proof we shall need some basic properties of $l$.

\begin{lemma}\label{lemma:lrr}
Suppose $P$ is a left-cancellative
semigroup with identity, $E$ is a product system over $P$,
$l:E\to\Bb(S(E))$ is the left regular representation,
$v,w\in E$, and $s\in P$. Then

\textup{(1)} $l(v)^*w$ is zero unless $p(w)\in p(v)P$.

\textup{(2)} If $p(w) = p(v)r$ for some $r\in P$, then
   $l(v)^*w \in E_r \subset S(E)$.

\textup{(3)} $\alpha^l_s(I)$ is the orthogonal projection onto
$\bigoplus_{t\in sP} E_t$.
\end{lemma}

\begin{proof} (1) Suppose $p(w) \notin p(v)P$. Then for any $u\in E$ we have
$p(w) \ne p(v)p(u) = p(vu)$,
so $\ip{l(v)^*w}u = \ip w{vu} = 0$.
Thus $l(v)^*w = 0$.

(2) If $u\in E_{p(v)}$ and $z \in E_r$,
then $l(v)^*(uz) = l(v)^*l(u)z = \ip uv z \in E_r$.
Since vectors of the form $uz$ have dense linear span in $E_{p(w)}$,
this gives (2).

(3) Let $\Bb$ be an orthonormal basis for $E_s$.
By (1) above,
\[
\alpha^l_s(I)w = \sum_{e\in\Bb} l(e)l(e)^*w = 0
\]
unless $p(w) \in sP$.
If $w = uz$ with $u\in\Bb$, $z\in E$, then
\[
\alpha^l_s(I)uz = \sum_{e\in\Bb} l(e)l(e)^*l(u)z
  = \sum_{e\in\Bb} \ip ue l(e)z = uz.
\]
Since vectors of this form are total in $\bigoplus_{t \in sP} E_t$,
this gives (3).
\end{proof}

\begin{proof}[Proof of Proposition~\ref{prop:lrr covariant}]
From Lemma~\ref{lemma:lrr}(3) we deduce that $\alpha^l_s(I)\alpha^l_t(I)$ is the
projection onto $\bigoplus\{E_r:r\in sP\cap tP\}$. But
\[
r\in sP\cap tP\Longleftrightarrow r\geq s\text{ and }r\geq t
\Longleftrightarrow r\geq s\vee t,
\]
so this is precisely the range of $\alpha^l_{s\vee t}(I)$.
\end{proof}

Since we shall be doing a lot of calculations with covariant representations, we
shall give some basic properties, and an alternative characterisation.

\begin{lemma}\label{lemma:handy}  Suppose $(G,P)$ is a quasi-lattice ordered group,
$E$ is a product system over $P$, $\phi$ is a representation of $E$
on $\Hh$, $u\in E$ and $s\in P$.

\textup{(1)} If $p(u) \le s$, then $\alpha^\phi_s(A)\phi(u) =
\phi(u)\alpha^\phi_{p(u)^{-1}s}(A)$ for any $A\in\bh$.

\textup{(2)} If $\phi$ is covariant, then
\[
\alpha^\phi_s(I)\phi(u) =
\begin{cases}
  \phi(u)\alpha^\phi_{p(u)^{-1}(p(u) \vee s)}(I)
  & \text{if $p(u) \vee s < \infty$,} \\
0 & \text{otherwise.}
\end{cases}
\]
\end{lemma}

\begin{proof} Suppose $p(u) \le s$.
Since $\alpha^\phi_{p(u)}(A)\phi(u) = \phi(u)A$ for each $A\in\bh$,
\[
\alpha^\phi_s(A)\phi(u) = \alpha^\phi_{p(u)}\bigl( \alpha^\phi_{p(u)^{-1}s}(A) \bigr) \phi(u)
  = \phi(u)\alpha^\phi_{p(u)^{-1}s}(A),
\]
giving (1).  If $\phi$ is covariant, then
$\alpha^\phi_s(I)\phi(u) = \alpha^\phi_s(I)\alpha^\phi_{p(u)}(I)\phi(u)$
is zero unless $p(u) \vee s < \infty$,
in which case
\begin{align*}
\alpha^\phi_s(I)\phi(u)
& = \alpha^\phi_s(I)\alpha^\phi_{p(u)}(I)\phi(u) \\
& = \alpha^\phi_{s \vee p(u)}(I) \phi(u) \\
& = \phi(u)\alpha^\phi_{p(u)^{-1}(p(u) \vee s)}(I),
\end{align*}
giving (2).
\end{proof}

\begin{prop}\label{prop:v*w}
Suppose $(G,P)$ is a quasi-lattice ordered group,
$E$ is a product system over $P$
and $\phi$ is a representation of $E$ on $\Hh$.

\textup{(1)} Suppose $v,w\in E$ satisfy $p(v) \vee p(w) < \infty$, and $\Bb$,
$\Cc$ are orthonormal bases for $E_{p(v)^{-1}(p(v) \vee p(w))}$
and $E_{p(w)^{-1}(p(v) \vee p(w))}$, respectively. Then the series
\[
\sum_{f\in \Bb, g\in\Cc} \ip{wg}{vf} \phi(f)\phi(g)^*,
\]
converges $\sigma$\ndash weakly to a bounded
operator on $\Hh$.

\textup{(2)} $\phi$ is covariant if and only if
for every $v,w\in E$
\begin{equation}\label{eq:v*w}
\phi(v)^*\phi(w) =
\begin{cases}
  \sum_{f,g} \ip{wg}{vf} \phi(f)\phi(g)^*
    & \text{if $p(v) \vee p(w) < \infty$} \\
  0 
    & \text{otherwise.}
\end{cases}
\end{equation}
\end{prop}

\begin{remark}\label{remark:norm convergent}
If $(G,P)$ is totally ordered, then either $f$ or $g$ disappears
from the sum in \eqref{eq:v*w}, and thus the series is norm convergent.
\end{remark}

\begin{proof}[Proof of Proposition~\ref{prop:v*w}]
(1) It does no harm to assume that $v$ and $w$ are unit vectors.
Then the series $\sum_f \rankone{vf}{vf}$ and $\sum_g \rankone {wg}{wg}$
converge strongly in the unit ball of $\Bb(E_{p(v) \vee p(w)})$,
and thus the series
\begin{equation}\label{eq:v*w2}
\sum_{f,g} (\rankone{vf}{vf})(\rankone{wg}{wg})
  = \sum_{f,g} \ip{wg}{vf} \rankone{vf}{wg}
\end{equation}
converges strongly to a bounded operator on $E_{p(v) \vee p(w)}$.
Since this convergence also occurs in the unit ball,
the series converges $\sigma$\ndash weakly.
Applying the isomorphism $\rho^\phi_{p(v) \vee p(w)}$
gives that the series $\sum_{f,g} \ip{wg}{vf} \phi(vf)\phi(wg)^*$
converges $\sigma$\ndash weakly, and multiplying on the left by $\phi(v)^*$
and on the right by $\phi(w)$ gives (1).

(2) If $\phi$ is covariant, then
\[
\phi(v)^*\phi(w)
= \phi(v)^*\alpha^\phi_{p(v)}(I)\alpha^\phi_{p(w)}(I)\phi(w)
\]
is zero unless $p(v) \vee p(w) < \infty$, in which case
\begin{align*}
\phi(v)^*\phi(w)
& = \phi(v)^*\alpha^\phi_{p(v) \vee p(w)}(I)\phi(w) \\
& = \alpha^\phi_{p(v)^{-1}(p(v) \vee p(w))}(I)
      \phi(v)^*\phi(w)
      \alpha^\phi_{p(w)^{-1}(p(v) \vee p(w))}(I) \\
& = \Big( \sum_f \phi(f)\phi(f)^* \Big)
      \phi(v)^*\phi(w)
      \Bigl( \sum_g \phi(g)\phi(g)^* \Bigr) \\
& = \sum_{f,g} \phi(f)\phi(vf)^* \phi(wg)\phi(g)^* \\
& = \sum_{f,g} \ip{wg}{vf} \phi(f)\phi(g)^*.
\end{align*}

Conversely, suppose \eqref{eq:v*w} holds for every $v,w\in E$.
Let $s,t\in P$.  Summing over $v,w$ in orthonormal bases for $E_s$
and $E_t$, respectively, we find that
\begin{align*}
\alpha^\phi_s(I)\alpha^\phi_t(I)
& = \Bigl( \sum_v \phi(v)\phi(v)^* \Bigr)
      \Bigl( \sum_w \phi(w)\phi(w)^* \Bigr) \\
& = \sum_{v,w} \phi(v)\phi(v)^*\phi(w)\phi(w)^*
\end{align*}
is zero unless $s \vee t < \infty$, in which case
\begin{align*}
\alpha^\phi_s(I)\alpha^\phi_t(I)
& = \sum_{v,w} \phi(v)\Bigl(
  \sum_{f,g} \ip{wg}{vf} \phi(f)\phi(g)^*
       \Bigr)\phi(w)^* \\
& = \sum_{v,w} \sum_{f,g}
  \ip{wg}{vf} \phi(vf)\phi(wg)^* \\
& = \sum_{v,w} \sum_{f,g}
 \phi(vf)\phi(vf)^*\phi(wg)\phi(wg)^* \\
& = \rho^\phi_{s \vee t}\Bigl( \sum_{v,w} \sum_{f,g}
 (\rankone{vf}{vf})(\rankone{wg}{wg}) \Bigr) \\
& = \rho^\phi_{s \vee t}(I) \\
& = \alpha^\phi_{s \vee t}(I),
\end{align*}
as required.
\end{proof}

\section{The system $(B_P, P, \tau, E)$}\label{section:system}

Suppose $(G,P)$ is a quasi-lattice ordered group
and $E$ is a product system over $P$.
For each $t\in P$ denote by $1_t$ the projection
in $\ell^\infty(P)$ defined by
\[
1_t(s) = 
\begin{cases}
  1 & \text{if $s \ge t$} \\
  0 & \text{otherwise.}
\end{cases}
\]
The product $1_s 1_t$ is $1_{s \vee t}$
if $s \vee t < \infty$ and $0$ otherwise; it follows that
$\Span\set{1_t: t\in P}$ is a \Star algebra, whose
closure is a \cstar subalgebra $B_P$ of $\ell^\infty(P)$.
The action  of $P$ by left translation on $\ell^\infty(P)$
restricts to an action $\tau$ of $P$ on $B_P$ such that
$\tau_s(1_t) = 1_{st}$ for $s,t\in P$.
We are interested in the twisted system $(B_P,P,\tau, E)$ because its
covariant representations are in one-to-one correspondence with the
covariant representations of $E$.

\begin{prop}\label{prop:pi phi}
Suppose $(G,P)$ is a quasi-lattice ordered group and $E$
is a product system over $P$.

\textup{(1)} If $(\pi, \phi)$ is a covariant representation of
$(B_P, P, \tau, E)$, then $\phi$ is a covariant representation of $E$ and
$\pi(1_s) = \alpha^\phi_s(I)$.

\textup{(2)} If $\phi$ is a covariant representation of $E$, then there
is a representation $\pi_\phi$  of $B_P$ such that
$\pi_\phi(1_s) = \alpha^\phi_s(I)$; moreover, $(\pi_\phi, \phi)$ is 
then a covariant representation of $(B_P, P, \tau, E)$.

\textup{(3)} $\pi_\phi$ is faithful iff
$\prod_{k=1}^n \bigl(I - \alpha^\phi_{s_k}(I)\bigr) \ne 0$
whenever $s_1, \dots, s_n\in P\setminus\set{e}$.
\end{prop}

\begin{proof} (1)  If $(\pi, \phi)$ is covariant,
then $\alpha^\phi_s(I) = \pi(\tau_s(1)) = \pi(1_s)$, so the
covariance of $\phi$ follows from the identity $1_s1_t=1_{s\vee t}$.

(2)  If $\phi$ is a covariant representation of $E$,
then by \cite[Proposition~1.3]{lacarae}
the map $1_s \mapsto \alpha^\phi_s(I)$ extends uniquely to a representation
$\pi_\phi$ of $B_P$.
Since $\pi_\phi(\tau_s(1_t)) = \pi_\phi(1_{st}) = \alpha^\phi_{st}(I)
          = \alpha^\phi_s(\alpha^\phi_t(I)) = \alpha^\phi_s(\pi_\phi(1_t))$,
$(\pi_\phi, \phi)$
is a covariant representation of $(B_P, P, \tau, E)$.

(3) By \cite[Proposition~1.3]{lacarae}, it suffices to show that
\begin{equation}\label{six}
\prod_{k=1}^n \bigl( \alpha^\phi_a(I) - \alpha^\phi_{z_k}(I) \bigr) \ne 0
\end{equation}
whenever $a, z_1, \dots, z_n \in P$ and $a < z_k$ for $k = 1,\dots, n$.
But $a<z_k$ means $z_k=as_k$ for some $s_k\in P\setminus \{e\}$, and
\[
\prod_{k=1}^n \bigl( \alpha^\phi_a(I) - \alpha^\phi_{z_k}(I) \bigr)
 = \alpha^\phi_a \Bigl( \prod_{k=1}^n \bigl( I - \alpha^\phi_{s_k}(I) 
\bigr) \Bigr),
\]
so the injectivity of $\alpha^\phi_a$ implies that
(\ref{six}) is equivalent to (3).
\end{proof}

\begin{cor}\label{cor:pi phi}
If $(G,P)$ is a quasi-lattice ordered group
and $E$ is a product system over $P$,
then the system $(B_P, P, \tau, E)$ has a covariant representation,
and  $i_{B_P}: B_P \to \bpp$
is faithful.
\end{cor}

\begin{proof} Since the left regular representation $l: E \to
\Bb(S(E))$ is covariant (Proposition~\ref{prop:lrr covariant}), the pair
$(\pi_l, l)$ is a covariant representation of $(B_P, P, \tau, E)$.
Lemma~\ref{lemma:lrr}(3) implies that the identity $\Omega$ of $E$,
viewed as an element of $E_e\subset S(E)$, is in the range of the
projection
$I - \alpha^l_s(I)$ whenever $s\in P \setminus\set{e}$,
and hence $\pi_l$ is faithful.
Since $\pi_l$ factors through $i_{B_P}$,
this in turn implies that $i_{B_P}$ is faithful.
\end{proof}

\begin{theorem}\label{theorem:subalgebra}
Suppose $(G,P)$ is a quasi-lattice ordered group and
$E$ is a product system over $P$. The \cstar subalgebra $A$ of $\bpp$
generated by the range of the canonical embedding $i_E$ is universal for
covariant representations of $E$, in the sense that:

\textup{(a)} there is a faithful unital representation $\sigma$  of $A$  on
Hilbert space such that $\sigma\circ i_E$ is a covariant
representation of
$E$, and

\textup{(b)} for every covariant representation $\phi$ of $E$ there is a
unital  representation $\pi$ of $A$ such that $\phi=\pi\circ i_E$.

\noindent If $E$ has finite-dimensional fibres, the algebra $A$ is all
of $\bpp$, and
\begin{equation}\label{eq:generators}
\bpp = \clsp\set{i_E(u)i_E(v)^*: u,v\in E}.
\end{equation}
\end{theorem}

\begin{remark}
Since the usual argument shows that there is at most one pair $(A,i_E)$
with these properties (see the proof of Proposition~\ref{prop:existence of cp}),
we can reasonably write $\cov(P,E)$ for $A:=C^*(i_E(E))$.
\end{remark}

\begin{remark}\label{remarks:C(E)}
As in \S\ref{section:crossed products}, one would expect and prefer to be able to replace
condition (a) by something like

(a$'$) for every representation $\sigma$ of $A$, 
$\sigma \circ i_E$ is a covariant representation of $E$.

\noindent If the sum in \eqref{eq:v*w} is always norm
convergent, then Proposition~\ref{prop:v*w} implies that conditions
(a) and (a$'$) are both equivalent to the following:
\[
i_E(v)^*i_E(w) =
\begin{cases}
  \sum_{f,g} \ip{wg}{vf} i_E(f)i_E(g)^*
    & \text{if $p(v) \vee p(w) < \infty$} \\
  0 
    & \text{otherwise,}
\end{cases}
\]
where the sum runs through orthonormal bases for
$E_{p(v)^{-1}(p(v) \vee p(w))}$
and $E_{p(w)^{-1}(p(v) \vee p(w))}$.
This is the case when the fibres of $E$ are finite-dimensional,
or when $(G,P)$ is totally ordered (Remark~\ref{remark:norm convergent}).
For this class of product systems $\cov(P,E)$ is indeed universal.
In a subsequent paper we will study a larger class of product systems
for which (a$'$) holds in $\cov(P,E)$.
\end{remark}

\begin{proof}[Proof of Theorem~\ref{theorem:subalgebra}] 
We can
represent $\bpp$ faithfully on a Hilbert space $\Kk$ in such a way
that $(i_{B_P},i_E)$ becomes a covariant representation of $(B_P, P,
\tau, E)$, and then $i_E$ is a covariant representation of $E$ by
Proposition~\ref{prop:pi phi}(1). If $\phi$ is a covariant
representation of $E$, then Proposition~\ref{prop:pi phi}(2) gives us
a covariant representation $(\pi_\phi,\phi)$ of $(B_P,P,\tau,E)$, and
hence a representation $\pi_\phi\times\phi$ of $\bpp$ such that
$(\pi_\phi\times\phi)\circ i_E=\phi$. Restricting
$\pi_\phi\times\phi$ to $A$ gives the required representation $\pi$.

Suppose now that $s\in P$ and that $\dim E_s < \infty$.
If $\Bb$ is an orthonormal basis for $E_s$, then
\[
i_{B_P}(1_s)
= i_{B_P}(\tau_s(1)) = \alpha^{i_E}_s(i_{B_P}(1))= \sum_{u\in\Bb}
i_E(u)i_E(u)^*
\]
belongs to $\cov(P,E)$. Thus if all the fibres of $E$ are
finite-dimensional we have $\cov(P,E)=\bpp$.

To establish \eqref{eq:generators}, it suffices to show that
$\Span\set{i_E(u)i_E(v)^*: u,v\in E}$
is closed under multiplication.
But by Proposition~\ref{prop:v*w}, each product
$i_E(u)i_E(v)^*i_E(w)i_E(z)^*$
is zero unless $p(v) \vee p(w) < \infty$, in which case
it is a {\em finite\/} sum of operators of the form
$i_E(uf)i_E(zg)^*$.
\end{proof}

\section{Faithful Representations}\label{section:faithfulness}

Our characterisation of faithful representations of $\bpp$ requires an 
amenability hypothesis, which we shall discuss shortly, and a spanning 
hypothesis, which says that
\begin{equation}\label{eq:spanning}
\bpp = \clsp\set{i_E(u)i_{B_P}(1_s)i_E(v)^*:
u,v\in E,\ s\in P}.
\end{equation}
This spanning hypothesis is automatically satisfied if $E$ has
finite-dimensional fibres (Theorem~\ref{theorem:subalgebra}),
or if $G$ is totally ordered (in which case we can simplify monomials
using $i_{B_P}(1_s)i_E(u) = i_E(u)$
or $i_E(u)i_{B_P}(1_{p(u)^{-1}s})$,
and the norm convergent expansion \eqref{eq:v*w}).

When the enveloping group $G$ of $P$ is abelian, the system
$(B_P,P,\tau,E)$ is {\em amenable\/} if averaging over the dual
action $\widehat\tau$ of $\widehat G$ gives a faithful expectation
onto the fixed-point algebra. In general we use the dual coaction
$\delta$ of
$G$ on $\bpp$ (Proposition~\ref{prop:coaction}), and the canonical trace
$\rho$ on $C^*(G)$ extending $f\mapsto f(e):\ell^1(G)\to\CC$. Then
$\Phi_\delta:=(\id\otimes\rho)\circ\delta$ is a positive linear map of
norm one of $B:=\bpp$ onto the fixed-point algebra
$B^\delta:=\{b\in B:\delta(b)=b\otimes 1\}$
(see \cite[2.3]{Ng} or \cite[Lemma~1.3]{Q1}).
A quick look at the characterisation of the coaction $\delta$ on
generators shows that
\[
\Phi_\delta(i_E(u)i_{B_P}(1_s)i_E(v)^*) =
\begin{cases}
  i_E(u)i_{B_P}(1_s)i_E(v)^* & \text{if $p(u) = p(v)$} \\
  0                          & \text{otherwise,}
\end{cases}
\]
and under the spanning hypothesis (\ref{eq:spanning}) this characterises
$\Phi_\delta$. (This implies, incidentally, that the expectation
$\Phi_\delta$ is independent of the choice of enveloping group $G$.)
We say the system is {\em amenable\/} if $\Phi_\delta$
is faithful in the sense that $\Phi_\delta(b^*b)=0$ implies $b=0$. The
argument of \cite[Lemma~6.5]{lacarae} shows that if the enveloping
group $G$ is amenable, then  $(B_P,P,\tau,E)$ is amenable
in our sense; in the next section we shall give further examples in
which
$P$ and $G$ are free products.

We can now state our main theorem.

\begin{theorem}\label{theorem:faithfulness of representations}
Suppose $(G,P)$ is a quasi-lattice ordered group,
$(B_P,P,\tau,E)$ is an amenable twisted system which satisfies the
spanning hypothesis \eqref{eq:spanning},
and $\phi$ is a covariant representation of $E$.
Then $\pi_\phi \times \phi$ is a faithful representation
of $\bpp$ if and only if
\begin{equation}\label{eq:piphifaithful3}
\prod_{k=1}^n \bigl(I - \alpha^\phi_{s_k}(I)\bigr) \ne 0
  \qquad\text{whenever $s_1, \dots, s_n\in P\setminus\set{e}$.}
\end{equation}
\end{theorem}

One direction is trivial: if $\pi_\phi\times\phi$ is faithful,
then by Corollary~\ref{cor:pi phi} so is
$\pi_\phi=(\pi_\phi\times\phi)\circ i_{B_P}$, and then
\eqref{eq:piphifaithful3} follows from Proposition~\ref{prop:pi phi}(3).
For the other direction, we follow the strategy of \cite[\S3]{lacarae}.
We show that for systems which satisfy the spanning hypothesis,
faithfulness of $\pi_\phi$ is sufficient to construct a spatial version
$\Phi_\phi$ of $\Phi_\delta$ such that
\[
\begin{CD} \bpp @>{\pi_\phi \times \phi}>>
  \pi_\phi \times \phi(\bpp)  \\ 
 @VV{\Phi_\delta}V       @VV{\Phi_\phi}V\\
 (\bpp)^\delta   @>{\pi_\phi \times \phi}>>
  \pi_\phi \times \phi((\bpp)^\delta)
\end{CD}
\]
commutes (Proposition~\ref{prop:existence of Phiphi}).
We also show that $\pi_\phi\times\phi$ is faithful on the
fixed-point algebra (Proposition~\ref{prop:faithfulness on fpa}), and
the amenability of the system completes the chain
\[
\begin{array}{rcl}
\pi_\phi\times\phi(b)=0&\Longrightarrow&
\Phi_\phi\big(\pi_\phi\times\phi(b^*b)\big)=0\\
&\Longleftrightarrow& \pi_\phi\times\phi(\Phi_\delta(b^*b))=0\\
&\Longleftrightarrow&\Phi_\delta(b^*b)=0\\
&\Longrightarrow& b=0.
\end{array}
\]

We begin by recalling some conventions from \cite[Lemma~1.4]{lacarae}.
Suppose  $F$ is a finite subset of $P$.
For each subset $A$ of $F$, define a projection $Q_A$ in $B_P$ by
\begin{equation}\label{eq:QA}
Q_A = 
\begin{cases}
  1_{\sigma A}\prod_{t\in F\setminus A} ( 1 - 1_t )
    & \text{if $\sigma A < \infty$} \\
  0 & \text{otherwise,}
\end{cases}
\end{equation}
with the convention that $\sigma\emptyset = e$. 

\begin{remark}\label{remark:QA}
It can be routinely verified  that $Q_A(s) = 1$ iff $A=\{t\in F:t\leq s\}$.
Thus $\set{Q_A: A\subset F}$
is a decomposition of the identity into mutually orthogonal projections,
and $Q_A$ is nonzero iff $A$ is an initial segment of $F$ in the sense
that $\sigma A < \infty$ and $A = \{t\in F: t \le \sigma A\}$.
In this case,
\begin{align*}
Q_A
& = \prod_{\set{t\in F: \sigma A < \sigma A \vee t < \infty}}
            \left( 1_{\sigma A} - 1_{\sigma A \vee t} \right) \nonumber \\
& = \tau_{\sigma A} \biggl(
 \prod_{\set{t\in F: \sigma A < \sigma A \vee t < \infty}}
            \bigl( 1 - 1_{\sigma A^{-1}(\sigma A \vee t)} \bigr)
                      \biggr)
\end{align*}
Thus if $\phi$ is a covariant representation of $E$
and $A$ is an initial segment of $F$,
\begin{equation}\label{eq:Q'A}
\pi_\phi(Q_A)
 = \alpha^\phi_{\sigma A}
     \biggl(\prod_{\set{t\in F: \sigma A < \sigma A \vee t < \infty}}
     \bigl( I - \alpha^\phi_{\sigma A^{-1}(\sigma A \vee t)}(I) \bigr) \biggr).
\end{equation}
\end{remark}

The following technical lemma will be used in the proofs of both 
Proposition~\ref{prop:faithfulness on fpa} and
Proposition~\ref{prop:existence of Phiphi}.

\begin{lemma}\label{lemma:QA}
Suppose $(G,P)$ is a quasi-lattice ordered group,
$E$ is a product system over $P$,
$\phi$ is a covariant representation of $E$,
$F$ is a finite subset of $P$,
$A$ is an initial segment of $F$,
$u,v\in E$ and $s\in P$.
Let $a = \sigma A$, so that $A = \set{t \in F: t \le a}$.

\textup{(1)} If $p(u) = p(v)$,
then the operator $\phi(u)\alpha^\phi_s(I)\phi(v)^*$
is in the commutant of $\pi_\phi(B_P)$.
In particular, it commutes with $\pi_\phi(Q_A)$.

\textup{(2)} If $p(u)s$, $p(v)s\in F$, then
\begin{multline*}
\pi_\phi(Q_A)\phi(u)\alpha^\phi_s(I)\phi(v)^*\pi_\phi(Q_A) \\
  = \left\{ \begin{array}{l}
      \pi_\phi(Q_A)
        \phi(u)\alpha^\phi_{p(u)^{-1}a}(I)\alpha^\phi_{p(v)^{-1}a}(I)\phi(v)^*
        \pi_\phi(Q_A) \\
      \phantom{0} \qquad\qquad\qquad\qquad \text{if $p(u)s\le a$ and 
$p(v)s
\le a$}
\\
      0 \qquad \text{otherwise.}
\end{array}
\right.
\end{multline*}
\end{lemma}

\begin{proof} (1) Suppose $p(u) = p(v)$;  it suffices to show that
$\phi(u)\alpha^\phi_s(I)\phi(v)^*$ commutes with
$\pi_\phi(1_t)$ for each $t\in P$.
If $p(u)s$ and $t$ have no common upper bound, then by Lemma~\ref{lemma:handy}
\begin{align*}
\pi_\phi(1_t)\phi(u)\alpha^\phi_s(I)\phi(v)^*
& = \alpha^\phi_t(I)\alpha^\phi_{p(u)s}(I)\phi(u)\phi(v)^* \\
& = 0 \\
& = \phi(u)\phi(v)^*\alpha^\phi_{p(v)s}(I)\alpha^\phi_t(I) \\
& = \phi(u)\alpha^\phi_s(I)\phi(v)^*\pi_\phi(1_t).
\end{align*}
Otherwise
\begin{align*}
\pi_\phi(1_t)\phi(u)\alpha^\phi_s(I)\phi(v)^*
& = \alpha^\phi_t(I) \phi(u)\alpha^\phi_s(I)\phi(v)^* \\
& = \phi(u)\alpha^\phi_{p(u)^{-1}(p(u) \vee t)}(I)\alpha^\phi_s(I)\phi(v)^* \\
& = \phi(u)\alpha^\phi_s(I)\alpha^\phi_{p(u)^{-1}(p(u) \vee t)}(I)\phi(v)^* \\
& = \phi(u)\alpha^\phi_s(I)\phi(v)^*\alpha^\phi_t(I).
\end{align*}

(2) Suppose $p(u)s, p(v)s \in F$.
The operator
\[
\pi_\phi(Q_A)\phi(u)\alpha^\phi_s(I)
 = \pi_\phi(Q_A)\alpha^\phi_a(I)\alpha^\phi_{p(u)s}(I)\phi(u)
\]
is zero unless $a \vee p(u)s < \infty$.
If $a < a \vee p(u)s < \infty$,
then
\[
\pi_\phi(Q_A) \le \pi_\phi(1_a - 1_{a \vee p(u)s})
  = \alpha^\phi_a(I) - \alpha^\phi_{a \vee p(u)s}(I),
\]
so that
\[
\pi_\phi(Q_A)\phi(u)\alpha^\phi_s(I)
= \pi_\phi(Q_A) \bigl( \alpha^\phi_a(I) - \alpha^\phi_{a \vee p(u)s}(I) \bigr)
               \alpha^\phi_{p(u)s}(I)\phi(u) = 0.
\]
Thus  $\pi_\phi(Q_A)\phi(u)\alpha^\phi_s(I)$
is zero unless $p(u)s \le a$, in which case
\begin{align*}
\pi_\phi(Q_A)\phi(u)\alpha^\phi_s(I)
& = \pi_\phi(Q_A)\alpha^\phi_a(I)\alpha^\phi_{p(u)s}(I)\phi(u) \\
& = \pi_\phi(Q_A)\alpha^\phi_a(I)\phi(u) \\
& = \pi_\phi(Q_A)\phi(u)\alpha^\phi_{p(u)^{-1}a}(I).
\end{align*}
Similarly, $\alpha^\phi_s(I)\phi(v)^*\pi_\phi(Q_A) = 0$
unless $p(v)s \le a$,
in which case it is equal to
$\alpha^\phi_{p(v)^{-1}a}\phi(v)^*\pi_\phi(Q_A)$.
Combining these results gives (2).
\end{proof}

Suppose  $s,t\in P$ and $s \le t$.
For each $A \in \Bb(E_s)$ denote by
$\beta_{t,s}(A)$ the unique operator on $E_t$ such that
\[
\beta_{t,s}(A)(uv) = (Au)v\ \text{ for } u\in E_s, \ v \in E_{s^{-1}t};
\]
if we use the multiplication to identify $E_s\otimes E_{s^{-1}t}$ with
$E_t$, then $\beta_{t,s}(A)$ is by definition $A\otimes I$. Each
$\beta_{t,s}$ is a faithful normal \Star homomorphism of $\Bb(E_s)$ into
$\Bb(E_t)$, and for $r \le s \le t$ we have
$\beta_{t,r} = \beta_{t,s} \circ \beta_{s,r}$.
If $u, v \in E_s$,
then $\beta_{t,s}(\rankone uv) = \sum_f \rankone{uf}{vf}$,
where $f$ ranges over an orthonormal basis for $E_{s^{-1}t}$.
If $\phi$ is a representation of $E$
and $\rho^\phi_t$ is the faithful normal \Star homomorphism
of Proposition~\ref{prop:rho}, we thus have
\begin{equation}\label{eq:rhobeta}
\begin{split}
\rho^\phi_t(\beta_{t,s}(\rankone uv))
& = \rho^\phi_t \Bigl( \sum_f \rankone{uf}{vf} \Bigr) \\
& = \sum_f \phi(uf)\phi(vf)^* \\
& = \phi(u)\alpha^\phi_{s^{-1}t}(I)\phi(v)^*.
\end{split}
\end{equation}

\begin{prop}\label{prop:faithfulness on fpa}
Suppose $(G,P)$ is a quasi-lattice ordered group,
$E$ is a product system over $P$ which satisfies the
spanning hypothesis \eqref{eq:spanning},
and $\phi$ is a covariant representation of $E$ which satisfies
\eqref{eq:piphifaithful3}.
Then the representation $\pi_\phi \times \phi$
of $\bpp$ is isometric on $(\bpp)^\delta$.
\end{prop}

\begin{proof} Let $X$ be a nonzero element in $\bpp$ of the form
\[
X = \sum_{(u,s,v)\in J} i_E(u)i_{B_P}(1_s)i_E(v)^*,
\]
where $J$ is a finite subset of
$\set{(u,s,v) \in E \times P \times E: p(u) = p(v)}$.
The spanning hypothesis \eqref{eq:spanning}
implies that elements such as $X$ are dense in
$(\bpp)^\delta$, so it suffices to show that
\begin{equation}\label{eq:piphiX}
\norm{\pi_\phi \times \phi(X)} = \norm X.
\end{equation}
Let $\sigma$ be a faithful representation of
$\bpp$ such that $(\sigma \circ i_{B_P}, \sigma \circ i_E)$
is a covariant representation of $(B_P, P, \tau, E)$.
By Proposition~\ref{prop:pi phi},
$i := \sigma \circ i_E$ is a covariant representation of $E$
and $\sigma \circ i_{B_P} = \pi_i$; in particular
$\pi_i(1_s) = \alpha^i_s(I)$ for each $s\in P$.

Let $F = \set{p(u)s: (u,s,v) \in J}$.
By Lemma~\ref{lemma:QA}, the operator $\pi_i \times i(X)$
commutes with each $\pi_i(Q_A)$.
Since these projections form a decomposition of the identity,
there is a subset $A \subseteq F$ such that
\[
\norm{\pi_i(Q_A)\pi_i \times i(X)} = \norm{\pi_i \times i(X)} = \norm X.
\]
Since $X \ne 0$, we have $\pi_i(Q_A)\not=0$.  From Remark~\ref{remark:QA} we
deduce that $a := \sigma A < \infty$ and $A$ is the initial
segment $\set{t\in F: t \le a}$.

Let $K:=\{(u,s,v)\in J:p(u)s\le a\}$, and define $T \in \Bb(E_a)$ by
\[
T = \sum_{(u,s,v)\in K} \beta_{a, p(u)} (\rankone uv).
\]
We claim that
\begin{equation}\label{eq:T}
  \norm{\pi_\phi \times \phi(X)} \ge \norm T = \norm X,
\end{equation}
from which \eqref{eq:piphiX} is immediate.
Suppose $\psi$ is a covariant representation of $E$;
we shall later take $\psi=i$ and $\psi=\phi$.
By Lemma~\ref{lemma:QA} and \eqref{eq:rhobeta},
\begin{align*}
\pi_\psi(Q_A)\pi_\psi \times \psi(X)
& = \pi_\psi(Q_A) \sum_{(u,s,v)\in J} \psi(u)\alpha^\psi_s(I)\psi(v)^* \\
& = \pi_\psi(Q_A) \sum_{(u,s,v)\in K}
                      \psi(u)\alpha^\psi_{p(u)^{-1}a}(I)\psi(v)^* \\
& = \pi_\psi(Q_A) \sum_{(u,s,v)\in K}
                      \rho^\psi_a\left( \beta_{a, p(u)}(\rankone uv) \right) \\
& = \pi_\psi(Q_A) \rho^\psi_a(T).
\end{align*}
From \eqref{eq:Q'A} we see that $\pi_\psi(Q_A)$
is in the range of $\alpha^\psi_a$.  If $\pi_\psi(Q_A)\not=0$,
Proposition~\ref{prop:rho} implies that
$S \mapsto \pi_\psi(Q_A)\rho^\psi_a(S)$ is a faithful representation
of $\Bb(E_a)$, so that
\[
\norm{\pi_\psi(Q_A)\pi_\psi \times \psi(X)}
  = \norm{\pi_\psi(Q_A)\rho^\psi_a(T)}
  = \norm T.
\]
We have already seen that $\pi_i(Q_A)\not=0$, and since
\eqref{eq:piphifaithful3} implies that $\pi_\phi$ is faithful, we
deduce that
 $\pi_\phi(Q_A)$ is nonzero.
Thus
\[
\norm{\pi_\phi \times \phi(X)}
  \ge \norm{\pi_\phi(Q_A)\pi_\phi \times \phi(X)}
  = \norm T
  = \norm{\pi_i(Q_A)\pi_i \times i(X)}
  = \norm X,
\]
so that \eqref{eq:T} holds as claimed.
\end{proof}

\begin{prop}\label{prop:existence of Phiphi}
Suppose $(G,P)$ is a quasi-lattice ordered group,
$E$ is a product system over $P$ which satisfies the
spanning hypothesis \eqref{eq:spanning},
and $\phi$ is a covariant representation of $E$ which satisfies
\eqref{eq:piphifaithful3}.
Let $\Delta = \set{(u,s,v) \in E \times P \times E: p(u) = p(v)}$.
Then there is a linear map $\Phi_\phi$ of norm one of
$\pi_\phi \times \phi(\bpp)$ onto $\pi_\phi \times \phi\bigl((\bpp)^\delta\bigr)$
such that, for each finite subset $J$ of $E \times P \times E$,
\[
\Phi_\phi\biggl(\sum_{(u,s,v) \in J}
\phi(u)\alpha^\phi_s(I)\phi(v)^*\biggr)
=
\sum_{(u,s,v) \in J \cap \Delta} \phi(u)\alpha^\phi_s(I)\phi(v)^*.
\]
\end{prop}

\begin{proof} Fix a finite subset $J$ of $E \times P \times E$ and let
\[
X = \sum_{(u,s,v)\in J} \phi(u)\alpha^\phi_s(I)\phi(v)^*,
\qquad
X_\Delta = \sum_{(u,s,v)\in J\cap\Delta} \phi(u)\alpha^\phi_s(I)\phi(v)^*.
\]
We will show that $\norm{X_\Delta} \le \norm X$, so that $\Phi_\phi$
is well-defined on finite sums such as $X$ and extends to a projection
of norm one on their closure, which by the spanning hypothesis
\eqref{eq:spanning} is all of $\pi_\phi \times \phi(\bpp)$.
Certainly we may assume that $X_\Delta \ne 0$.
Let 
\[
F = \set{p(u)s: (u,s,v)\in J} \cup \set{p(v)s: (u,s,v) \in
J}.
\]
By Lemma~\ref{lemma:QA}, $X_\Delta$ commutes with each
$\pi_\phi(Q_A)$, and  since the $Q_A$ form a
decomposition of the identity, there exists $A\subseteq F$ such that
\[
\norm{X_\Delta} = \norm{\pi_\phi(Q_A) X_\Delta}.
\]
Because $X_\Delta \ne 0$, we have
$\pi_\phi(Q_A)\not=0$, so by Remark~\ref{remark:QA}, we have $a := \sigma A < \infty$
and $A=\set{t\in F: t \le a}$.
It follows from (\ref{eq:Q'A}) that $\pi_\phi(Q_A)$
is in the range of $\alpha^\phi_a$.

Define $T\in\Bb(E_a)$ by
\[
T = \sum_{\{(u,s,v)\in J \cap \Delta :p(u)s \le a\}}
         \beta_{a,p(u)}(\rankone uv).
\]
Exactly as in the proof of Proposition~\ref{prop:faithfulness on fpa}
we have $\pi_\phi(Q_A) X_\Delta = \pi_\phi(Q_A) \rho^\phi_a(T)$,
so that by Proposition~\ref{prop:rho} we have
\[
\norm{X_\Delta}
 = \norm{\pi_\phi(Q_A) X_\Delta}
 = \norm{\pi_\phi(Q_A) \rho^\phi_a(T)}
 = \norm T.
\]
We will construct another nonzero projection $Q$ in the range of $\alpha^\phi_a$
with the property that $QXQ = Q\rho^\phi_a(T)$.  This will complete
the proof, since from this another application of Proposition~\ref{prop:rho}
gives 
\[
\norm{X_\Delta} = \norm T = \norm{Q\rho^\phi_a(T)} = 
\norm{QXQ} \le \norm X.
\]

For each $b,c\in A$ such that $b \ne c$ and
$b^{-1}a \vee c^{-1}a < \infty$,
define $d_{b,c}\in P$ as in \cite[Lemma~3.2]{lacarae}:
\[
d_{b,c} =
\begin{cases}
  (b^{-1}a)^{-1}(b^{-1}a \vee c^{-1}a)
    & \text{if $b^{-1}a < b^{-1}a \vee c^{-1}a$} \\
  (c^{-1}a)^{-1}(b^{-1}a \vee c^{-1}a)
    & \text{otherwise,} \\
\end{cases}
\]
noting in particular that $d_{b,c}$ is never the identity in $P$.
Define
\[
R' = \prod_{\substack{b\ne c\in A\\ b^{-1}a \vee c^{-1}a <
\infty}}
\bigl( I - \alpha^\phi_{d_{b,c}}(I) \bigr),
\]
\[
Q' = \prod_{\substack{t\in F \\ a < a \vee t < \infty}}
           \bigl( I - \alpha^\phi_{a^{-1}(a \vee t)}(I) \bigr)
     \prod_{\substack{b\ne c\in A\\ b^{-1}a \vee c^{-1}a < \infty}}
           \bigl(I - \alpha^\phi_{d_{b,c}}(I) \bigr),
\]
$R = \alpha^\phi_a(R')$ and $Q = \alpha^\phi_a(Q')$, so that $Q = \pi_\phi(Q_A)R$.
By condition \eqref{eq:piphifaithful3}, $Q' \ne 0$, and thus
$Q$ is a nonzero projection in the range of $\alpha^\phi_a$.
We claim that $Q$ is the desired projection satisfying
$QXQ = Q\rho^\phi_a(T)$.

To begin with, because $Q \le \pi_\phi(Q_A)$,
we can use  Lemma~\ref{lemma:QA} to rewrite
\begin{equation}\label{eq:QXQ}
 QXQ= Q \biggl( \sum_{\substack{(u,s,v) \in J \\ p(u)s, p(v)s \le a}}
  \phi(u)\alpha^\phi_{p(u)^{-1}a}(I)\alpha^\phi_{p(v)^{-1}a}(I)\phi(v)^* \biggr) Q.
\end{equation}
Now suppose $(u,s,v) \in J$, $p(u)s\le a$, $p(v)s \le a$ and $p(u) \ne
p(v)$. If $p(u)^{-1}a$ and $p(v)^{-1}a$ have no common upper bound,
then the corresponding term in the above sum is zero.
On the other hand, if $p(u)^{-1}a \vee p(v)^{-1}a < \infty$,
then
\[
(p(u)s)^{-1}a \vee (p(v)s)^{-1}a
 = s^{-1}(p(u)^{-1}a \vee p(v)^{-1}a)<\infty.
\]
Let $d = d_{p(u)s, p(v)s}$.
The previous equation shows that $d$ is either
\[
(p(u)^{-1}a)^{-1}(p(u)^{-1}a \vee p(v)^{-1}a)\ 
\text{ or }\ 
(p(v)^{-1}a)^{-1}(p(u)^{-1}a \vee p(v)^{-1}a).
\]
Then $R \le \alpha^\phi_a(I) - \alpha^\phi_{ad}(I)$, and
\begin{multline*}
\bigl( \alpha^\phi_a(I) - \alpha^\phi_{ad}(I) \bigr)
  \phi(u)\alpha^\phi_{p(u)^{-1}a}(I)\alpha^\phi_{p(v)^{-1}a}(I)\phi(v)^*
  \bigl( \alpha^\phi_a(I) - \alpha^\phi_{ad}(I) \bigr) \\
\begin{split}
  & = \phi(u)\bigl( \alpha^\phi_{p(u)^{-1}a}(I) - \alpha^\phi_{p(u)^{-1}ad}(I) \bigr)
      \alpha^\phi_{p(u)^{-1}a}(I) \\
  & \qquad\qquad\qquad \alpha^\phi_{p(v)^{-1}a}(I)
    \bigl( \alpha^\phi_{p(v)^{-1}a}(I) - \alpha^\phi_{p(v)^{-1}ad}(I) \bigr)\phi(v)^* \\
  & = 0,
\end{split} 
\end{multline*}
since either $p(u)^{-1}ad$ or $p(v)^{-1}ad$ is equal to $p(u)^{-1}a \vee p(v)^{-1}a$.
This shows that
\[
R\phi(u)\alpha^\phi_{p(u)^{-1}a}(I)\alpha^\phi_{p(v)^{-1}a}(I)\phi(v)^*R = 0
\]
for each
$(u,s,v) \in J$ satisfying $p(u)s\le a$, $p(v)s \le a$ and $p(u) \ne
p(v)$. Equation \eqref{eq:QXQ} now simplifies to
\[
QXQ
 = Q\biggl( \sum_{\substack{(u,s,v) \in J \\ p(u)s = p(v)s \le a}}
  \phi(u)\alpha^\phi_{p(u)^{-1}a}(I)\phi(v)^* \biggr)Q
= Q\rho^\phi_a(T),
\]
so this $Q$ will suffice.
\end{proof}

\begin{examples}\label{examples:faithfulness of representations}
(1) Applying Theorem~\ref{theorem:faithfulness of representations}
to the trivial product system $P\times \CC$ gives
\cite[Theorem~3.7]{lacarae}.
More generally, if $\mu$ is a multiplier on $P$, applying it
to $(P\times\CC)^\mu$  gives a characterisation of  the faithful
representations of the universal \cstar algebra for covariant
$\mu$\ndash representations of $(G,P)$.

(2) If $E$ is a product system over $\NN$ with $\dim E_1 = n <
\infty$, then $\cov(P,E)$ is the Toeplitz-Cuntz algebra $\Tt\Oo_n$. 
In this case, our Theorem~\ref{theorem:faithfulness of representations}
reduces to Cuntz's Theorem: the representation of $\Tt\Oo_n$
corresponding to a Toeplitz-Cuntz family
$\set{V_1, V_2, \dots, V_n}$ is faithful iff $\sum V_kV_k^* < I$.  

Since $(\ZZ, \NN)$ is totally ordered, the theorem still applies when
$\dim E_1 = \aleph_0$, and states that the representation of
$B_{\NN} \cross_{\tau, E} \NN$ corresponding to an infinite Toeplitz-Cuntz
family $\set{V_1, V_2, \dotsc}$ is faithful iff $\sum V_kV_k^* < I$;
this implies  in particular that $B_{\NN} \cross_{\tau, E} \NN$ is not
simple. Since $\cov(\NN,E)$ is isomorphic to the simple \cstar
algebra $\Oo_\infty$, it is a proper subalgebra of
$B_{\NN} \cross_{\tau, E} \NN$. The system considered in Example~\ref{ex:not
universal} is exactly $(B_{\NN}, \NN, \tau, E)$.

(3) Consider the lexicographic product system
$E$ over $\NN\oplus\NN$ determined by the homomorphism
$d: (m,n) \in \NN\oplus\NN \mapsto 2^m3^n \in \NN^*$.
A representation $\phi$ of $E$ is determined by
a pair of Toeplitz-Cuntz families
$\set{U_1, U_2}$, $\set{V_1, V_2, V_3}$ satisfying (\ref{commrel}),
and from Proposition~\ref{prop:v*w} it is easy to see
that $\phi$ is covariant iff
\begin{align*}
U_1^*V_1 & = V_1U_1^* + V_2U_2^*, & U_2^*V_1 & = 0, \\
U_1^*V_2 & = V_3U_1^*,            & U_2^*V_2 & = V_1U_2^*, \\
U_1^*V_3 & = 0,                   & U_2^*V_3 & = V_2U_1^* + V_3U_2^*.
\end{align*}
Thus $\cov(\NN^2,E)$ is universal  for pairs
of Toeplitz-Cuntz families satisfying all of these relations,
and Theorem~\ref{theorem:faithfulness of representations}
implies that $\set{U_1, U_2}$ and $\set{V_1, V_2, V_3}$
generate a faithful representation of $\cov(\NN^2,E)$ iff
\[
(I - U_1U_1^* - U_2U_2^*)(I - V_1V_1^* - V_2V_2^* - V_3V_3^*) \ne 0.
\]
\end{examples}

We conclude this section by showing that, under the spanning hypothesis
\eqref{eq:spanning}, amenability of $E$ (strictly speaking, amenability of
$(B_P,P,\tau,E)$) is equivalent to faithfulness of the left regular representation.

\begin{prop}\label{prop:Phil faithful}
Suppose $(G,P)$ is a quasi-lattice ordered group and $E$
is a product system over $P$.
Let $l: E \to \Bb(S(E))$ be the left regular representation of $E$ and let
\[
X = l(u_1)\alpha^l_{s_1}(I)l(v_1)^* \dotsm l(u_n)\alpha^l_{s_n}(I)l(v_n)^*.
\]
Then the map
\begin{equation}\label{eq:Phil}
X \mapsto
\begin{cases}
  X & \text{if $p(u_1)p(v_1)^{-1}\dotsm p(u_n)p(v_n)^{-1} = e$} \\
  0 & \text{otherwise}
\end{cases}
\end{equation}
extends to a projection $\Phi_l$ of norm one on
$\pi_l \times l(\bpp)$ which is faithful on positive elements.
\end{prop}

\begin{proof} For each $s\in P$ let $Q_s$ be the orthogonal projection
of $S(E)$ onto $E_s$.
Since the $Q_s$'s are mutually orthogonal, the formula
\[
\Phi_l(T) = \sum_{s\in P} Q_s T Q_s,\qquad T\in\Bb(S(E)),
\]
defines a completely positive projection of norm one on $\Bb(S(E))$
which is faithful on positive operators.
We claim that the restriction of $\Phi_l$
to $\pi_l \times l(\bpp)$
satisfies \eqref{eq:Phil}.

Let $r = p(u_1)p(v_1)^{-1} \dotsm p(u_n)p(v_n)^{-1}$.
For each $s \in P$, Lemma~\ref{lemma:lrr} implies that
$X$ is zero on $E_s$
unless $rs \in P$, in which case $X$ maps $E_s$ into $E_{rs}$.
Thus if $r \ne e$, $Q_sXQ_s=0$ for every $s\in P$,
and $\Phi_l(X) = 0$.  If on the other hand
$r = e$, then $Q_sXQ_s = XQ_s$ for each $s\in P$, and
\[
\Phi_l(X) = \sum_{s\in P} Q_sXQ_s = X\sum_{s\in P} Q_s = X.
\]
\end{proof}

\begin{cor}\label{cor:Phil faithful}
Suppose $(G,P)$ is a quasi-lattice ordered group and $E$ is
a product system over $P$ satisfying \eqref{eq:spanning}.
Then $E$ is amenable if and only if $\pi_l \times l$ is faithful.
\end{cor}

\begin{proof} Suppose $\pi_l \times l$ is faithful.
By Proposition~\ref{prop:Phil faithful},
$(\pi_l \times l) \circ \Phi_\delta = \Phi_l \circ (\pi_l \times l)$ is faithful
on positive elements, hence so is $\Phi_\delta$; that is, $E$ is amenable.
If \eqref{eq:spanning} is satisfied,
then Proposition~\ref{prop:faithfulness on fpa} implies that
$\pi_l \times l$ is faithful on $(\bpp)^\delta$.
If in addition $E$ is amenable, then
$\Phi_l \circ (\pi_l \times l) = (\pi_l \times l) \circ \Phi_\delta$
is faithful on positive elements, from which faithfulness of
$\pi_l \times l$ follows.
\end{proof}

\section{Amenability}\label{section:amenability}

Suppose $(G,P)$ is a quasi-lattice ordered group
and $E$ is a product system over $P$.
In this section we give conditions which ensure that $E$ is amenable; these
conditions also ensure that
$E$ satisfies the spanning condition \eqref{eq:spanning},
so that  Theorem~\ref{theorem:faithfulness of representations} applies.
Our argument follows those of \cite[\S4]{lacarae} and
\cite[Proposition~2.10]{dinhjfa}.

\begin{theorem}\label{theorem:amenability}
Suppose $\theta: (G, P) \to (\Gg, \Pp)$ is a homomorphism of quasi-lattice
ordered groups such that, whenever $s\vee t<\infty$,
\begin{equation}\label{eq:theta}
  \theta(s \vee t)=  \theta(s) \vee \theta(t)\ \text{ and }\ 
  \theta(s)=  \theta(t)\Longrightarrow s = t,
\end{equation}
and suppose that $\Gg$ is amenable.
If $E$ is a product system over $P$ which satisfies 
\begin{multline}\label{eq:v*w to fg*}
  i_E(v)^*i_E(w) \in \clsp\{i_E(f)i_E(g)^*: \\
  f\in E_{p(v)^{-1}(p(v) \vee p(w))}, g\in E_{p(w)^{-1}(p(v) \vee p(w))}\},
\end{multline}
then $E$ is amenable and the spanning hypothesis \eqref{eq:spanning} holds.
\end{theorem}

\begin{remark}
(1) As in \cite[Proposition~4.3]{lacarae},
the main example of such a map $\theta$ will be the canonical homomorphism of
a free product of quasi-lattice ordered groups
onto the corresponding direct sum. However, we could also take $\theta$ to be
the length function on the free group $\field{F}_n$ (the homomorphism into
$\ZZ$ which takes each generator to $1$), and this example gives a good feel for
both our  constructions and those of \cite[\S4]{lacarae}.
\end{remark}

\begin{proof}[Proof of Theorem~\ref{theorem:amenability}] 
The homomorphism $\theta:G \to \Gg$ induces a coaction
$\delta_\theta = (\id \otimes \theta) \circ \delta$ of $\Gg$ on $\bpp$,
and hence a conditional expectation $\Phi_{\delta_\theta}$
of $\bpp$ onto the fixed-point algebra $(\bpp)^{\delta_\theta}$,
such that 
\[
\Phi_{\delta_\theta}(i_E(u)i_{B_P}(1)i_E(v)^*) =
\begin{cases}
  i_E(u)i_{B_P}(1)i_E(v)^*& \text{if $\theta(p(u)) = \theta(p(v))$} \\
  0 & \text{otherwise.}
\end{cases}
\]
Since $\Gg$ is amenable, $\Phi_{\delta_\theta}$ is faithful on positive elements.
We can  recover the original expectation $\Phi_\delta$
by first applying $\Phi_{\delta_\theta}$,
and then killing the terms  with
$p(u)p(v)^{-1}\in\ker \theta\setminus\{e\}$,
which can be accomplished spatially by
representing $(\bpp)^{\delta_\theta}$
using the regular representation $\pi_l \times l$, and compressing
to the diagonal via the expectation $\Phi_l$
of Proposition~\ref{prop:Phil faithful}.
Since $\Phi_l$ is faithful on positive operators,
this last step is faithful whenever $\pi_l \times l$
is faithful on $(\bpp)^{\delta_\theta}$.

It therefore suffices to show that $\pi_l \times l$
is faithful on $(\bpp)^{\delta_\theta}$.
Let $\sigma$ be a faithful representation of
$\bpp$ such that $(\sigma \circ i_{B_P}, \sigma \circ i_E)$
is a covariant representation of $(B_P, P, \tau, E)$.
By Proposition~\ref{prop:pi phi},
$i = \sigma \circ i_E$ is a covariant representation of $E$
and $\sigma \circ i_{B_P} = \pi_i$; in particular, we have
$\pi_i(1_s) = \alpha^i_s(I)$ for each $s\in P$.

Suppose $S$ is a subset of $\Pp$ for which $q \vee r \in S$
whenever $q,r\in S$ and $q \vee r < \infty$.
We claim that
\[
\Uu_S = \clsp\set{i_E(u)i_{B_P}(1_s)i_E(v)^*: \theta(p(u)s) = \theta(p(v)s) \in S}
\]
is a \cstar subalgebra of $\bpp$.
For this, suppose that $u,v,w,z \in E$ and $s,t\in P$ are such that
$\theta(p(u)s) = \theta(p(v)s) \in S$
and $\theta(p(w)t) = \theta(p(z)t) \in S$. Then by Lemma~\ref{lemma:handy},
\[
i(u)\alpha^i_s(I)i(v)^*i(w)\alpha^i_t(I)i(z)^*
 = i(u)i(v)^*\alpha^i_{p(v)s}(I)\alpha^i_{p(w)t}(I)i(w)i(z)^*.
\]
By Proposition~\ref{prop:v*w}, this operator is
is zero unless $p(v)s \vee p(w)t < \infty$, in which case
by \eqref{eq:v*w to fg*} it can be approximated in norm
by a finite sum of operators of the form
\begin{equation}\label{summand}
i(u)\alpha^i_s(I)i(f)i(g)^*\alpha^i_t(I)i(z)^*,
\end{equation}
where $p(f) = p(v)^{-1}(p(v) \vee p(w))$
and $p(g) = p(w)^{-1}(p(v) \vee p(w))$.
Again using Lemma~\ref{lemma:handy}, each operator (\ref{summand}) can be
rewritten as
\begin{equation}\label{eq:am1}
i(uf)\alpha^i_{p(f)^{-1}(p(f) \vee s)}(I)\alpha^i_{p(g)^{-1}(p(g) \vee
t)}(I)i(zg)^*.
\end{equation}
Now
\begin{align*}
p(f)^{-1}(p(f) \vee s)
& = (p(v) \vee p(w))^{-1} p(v)(p(f) \vee s) \\
& = (p(v) \vee p(w))^{-1} (p(v)p(f) \vee p(v)s) \\
& = (p(v) \vee p(w))^{-1} (p(v) \vee p(w) \vee p(v)s) \\
& = (p(v) \vee p(w))^{-1} (p(v)s \vee p(w)),
\end{align*}
and similarly
$p(g)^{-1}(p(g) \vee t) = (p(v) \vee p(w))^{-1} (p(v) \vee p(w)t)$.
Thus
\begin{multline*}
p(f)^{-1}(p(f) \vee s) \vee p(g)^{-1}(p(g) \vee t) \\
\begin{split}
& = (p(v) \vee p(w))^{-1} (p(v)s \vee p(w))
      \vee (p(v) \vee p(w))^{-1} (p(v) \vee p(w)t) \\
& = (p(v) \vee p(w))^{-1} (p(v)s \vee p(w) \vee p(v) \vee p(w)t) \\
& = (p(v) \vee p(w))^{-1} (p(v)s \vee p(w)t).
\end{split}
\end{multline*}
Using this to simplify \eqref{eq:am1}, we see that
$i(u)\alpha^i_s(I)i(v)^*i(w)\alpha^i_t(I)i(z)^*$
can be approximated in norm by a finite sum of operators
of the form
\[
i(uf)\alpha^i_{(p(v) \vee p(w))^{-1} (p(v)s \vee p(w)t)}(I)i(zg)^*.
\]
Now
\begin{align*}
\theta(p(uf)(p(v) \vee p(w))^{-1} (p(v)s \vee p(w)t))
& = \theta(p(u)p(v)^{-1} (p(v)s \vee p(w)t)) \\
& = \theta(p(v)s \vee p(w)t),
\end{align*}
and similarly 
$\theta(p(zg)(p(v) \vee p(w))^{-1} (p(v)s \vee p(w)t)) = \theta(p(v)s \vee p(w)t)$.
Since\[
\theta(p(v)s \vee p(w)t) = \theta(p(v)s) \vee \theta(p(w)t) \in S,
\]
this shows that
$i(u)\alpha^i_s(I)i(v)^*i(w)\alpha^i_t(I)i(z)^*$
is an element of $\sigma(\Uu_S)$, and hence that $\Uu_S$ is closed under
multiplication. This proves that $\Uu_S$ is a \cstar algebra.

Minor revisions of the above argument show that
$\clsp\set{i_E(u)i_{B_P}(1_s)i_E(v)^*: u,v\in E, s\in P}$
is a \cstar algebra, so that \eqref{eq:spanning} holds.
Applying $\Phi_{\delta_\theta}$ to both sides of \eqref{eq:spanning} gives
\[
(\bpp)^{\delta_\theta}
  = \Uu_{\Pp} := \clsp\set{i_E(u)i_{B_P}(1_s)i_E(v)^*: \theta(p(u)) =
\theta(p(v))}.
\]
Let $\Ff$ be the set of all finite subsets of $\Pp$
which are closed under $\vee$.
As in \cite[Lemma~4.1]{lacarae}, $\Ff$ is directed under set inclusion,
so that
\[
(\bpp)^{\delta_\theta} = \overline{\bigcup_{F\in\Ff} \Uu_F}.
\]
By \cite[Lemma~1.3]{alnr}, to prove that $\pi_l \times l$ is faithful
on $(\bpp)^{\delta_\theta}$ it is enough to prove it is faithful
 on each of the subalgebras $\Uu_F$.
We shall accomplish this by inducting on $\abs F$.

First suppose $F = \set{r}$ for some $r\in\Pp$, and write $\Uu_r$
for $\Uu_{\set{r}}$.
Let $\phi$ be a covariant representation of $E$,
and suppose that $x$ and $y$ are unit vectors in $E$ such that
$p(x) \ne p(y)$ and $\theta(p(x)) = \theta(p(y)) = r$.
Since $\theta$ satisfies \eqref{eq:theta} we must have
$p(x) \vee p(y) = \infty$, and thus by Proposition~\ref{prop:v*w}
the isometries $\phi(x)$ and $\phi(y)$ have orthogonal ranges.
Hence $\phi$ extends to a bounded linear map on $F_r := \bigoplus_{t\in\theta^{-1}(r)}
E_t$, and the following analogues of Propositions~\ref{prop:alphaphi}
and \ref{prop:rho} hold:
$\alpha^\phi_r := \sum_{t\in\theta^{-1}(r)} \alpha^\phi_t$
defines a normal \Star endomorphism of $\Bb(\Hh_\phi)$,
and there is a faithful normal \Star representation $\rho^\phi_r$
of $\Bb(F_r)$ such that $\rho^\phi_r(\rankone xy) = \phi(x)\phi(y)^*$ for
$x,y\in F_r$.  Moreover, if $Q$ is a nonzero projection on $\Hh_\phi$,
then $T \mapsto \alpha^\phi_r(Q)\rho^\phi_r(T)$ is also a faithful
representation of $\Bb(F_r)$.  There should be no confusion caused
by our abuse of notation; just take note of whether the subscript
is an element of $P$ or $\Pp$.  A word of caution, however:
although $t \mapsto \alpha^\phi_t$ is a semigroup homomorphism,
in general  the map $r\in \Pp \mapsto \alpha^\phi_r$ is not:
the bundle $\{F_r:r\in\Pp\}$ is not a product system in the multiplication
inherited from $E$.

Suppose that $J$ is a finite subset of
$\set{(u,s,v) \in E \times P \times E: \theta(p(u)s) = \theta(p(v)s) = r}$,
and let
\[
X = \sum_{(u,s,v) \in J} i_E(u)i_{B_P}(1_s)i_E(v)^*;
\]
to prove $\pi_l \times l$ faithful on $\Uu_r$
we will show that $\norm{\pi_l \times l(X)} = \norm X$.
Define $T \in \Bb(F_r)$ by
\[
T = \sum_{(u,s,v) \in J} \sum_f \rankone{uf}{vf},
\]
where $f$ ranges over an orthonormal basis for $E_s$.
It is routine to check that
\[
\rho^l_r(T) = \sum_{(u,s,v) \in J} l(u)\pi_l(1_s)l(v)^* = \pi_l \times l(X),
\]
and similarly $\rho^i_r(T) = \pi_i \times i(X) = \sigma(X)$.
Since $\rho^l_r$, $\rho^i_r$ and $\sigma$ are isometric,
\[
\norm{\pi_l \times l(X)} = \norm{\rho^l_r(T)} = \norm T
  = \norm{\rho^i_r(T)} = \norm{\sigma(X)} = \norm X.
\]

For the inductive step,
suppose $F \in \Ff$ and $\pi_l \times l$ is faithful on $\Uu_{F'}$
whenever $F' \in \Ff$ and $\abs{F'} < \abs F$; we aim to prove that
$\pi_l \times l$ is faithful on $\Uu_F$. 
Since $F$ is finite it has a minimal element;
that is, there exists $r_0 \in F$ such that
$r_0 < r_0 \vee r$ for each $r\in F \setminus\set{r_0}$.
Notice that if $u,v,w \in E$ and $s\in P$ are such that
$\theta(p(u)s) = \theta(p(v)s) \in F$ and $\theta(p(w)) = r_0$,
then by Lemma~\ref{lemma:lrr} the vector
$\pi_l \times l(i_E(u)i_{B_P}(1_s)i_E(v)^*)w = l(u)\alpha^l_s(I)l(v)^*w$
is nonzero only when $p(v)s \le p(w)$.
Since this in turn implies that $\theta(p(v)s) \le r_0$,
the minimality of $r_0$ forces $\theta(p(v)s) = r_0$.
Thus if we let $P_{r_0}$ denote the orthogonal projection
of $S(E)$ onto $F_{r_0}$,
then $\pi_l \times l(\Uu_r)P_{r_0} = \set{0}$ for each
$r \in F \setminus\set{r_0}$.

On the other hand, we have already demonstrated that
$\pi_l \times l$ maps $\Uu_{r_0}$ into
the range of $\rho^l_{r_0}$, and an easy calculation
shows that $P_{r_0} = \alpha^l_{r_0}(Q_e)$,
where $Q_e$ is the orthogonal projection onto $E_e$.
Since $T \mapsto \alpha^l_{r_0}(Q)\rho^l_{r_0}(T)$
is a faithful normal \Star homomorphism,
so is the map $X \in \Uu_{r_0} \mapsto \pi_l \times l(X)P_{r_0}$.

Now suppose $Y \in \Uu_F$ and $\pi_l \times l(Y) = 0$.
We will show that $Y \in \Uu_{F \setminus \set{r_0}}$,
from which the inductive hypothesis implies that $Y = 0$.
Let $(Y_n)$ be a sequence in
\[
\Span\set{i_E(u)i_{B_P}(1_s)i_E(v)^*: \theta(p(u)s) = \theta(p(v)s) \in F}
\]
which converges in norm to $Y$, and express each $Y_n$ as a sum
$\sum_{r\in F} Y_{n,r}$, where $Y_{n,r} \in \Uu_r$.
For each $n$,
\[
\norm{\pi_l \times l(Y_n) P_{r_0}}
 = \norm{\pi_l \times l(Y_{n,r_0}) P_{r_0}}
 = \norm{Y_{n,r_0}},
\]
and consequently $Y_{n, r_0} \ra 0$.
Thus $Y_n - Y_{n,r_0} \ra Y$, which shows that
$Y \in \Uu_{F \setminus\set{r_0}}$, as claimed.
\end{proof}

\begin{cor}\label{cor:amenability}
Suppose $(G^\lambda, P^\lambda)$ is a quasi-lattice ordered
group with $G^\lambda$ am\-en\-able for each $\lambda$ belonging
to some index set $\Lambda$.
Then any product system over $* P^\lambda$
which satisfies \eqref{eq:v*w to fg*} is amenable.
In particular, any product system over $* P^\lambda$
which has only finite-dimensional fibres is amenable.
\end{cor}

\begin{proof} The group $\bigoplus G^\lambda$ is amenable,
and by \cite[Proposition~4.3]{lacarae} the canonical map
$\theta: * G^\lambda \to \bigoplus G^\lambda$ satisfies
\eqref{eq:theta}. It follows from Proposition~\ref{prop:v*w} that any system
with finite-dimensional fibres will satisfy \eqref{eq:v*w to fg*}.
\end{proof}

\end{document}